\documentclass[12pt,a4paper,twoside]{amsart}
\usepackage{amssymb,latexsym,amsmath,graphics}
\newcommand{\nnntext}{}
\newcommand{\mntext}{}
\newcommand{\ntext}{}

\newtheorem {theorem}{Theorem}
\newtheorem {lemma}[theorem]{Lemma}
\newtheorem {proposition}[theorem]{Proposition}
\newtheorem {definition}[theorem]{Definition}
\newtheorem {corollary}[theorem]{Corollary}
\newtheorem{remark}[theorem]{Remark}

\newcommand{\Meas}{\Theta}

\newcommand{\T}{{\mathbb T}}
\newcommand{\R}{{\mathbb R}}
\newcommand{\N}{{\mathbb N}}
\newcommand{\Z}{{\mathbb Z}}
\newcommand{\C}{{\mathbb C}}
\newcommand{\B}{{\mathcal B}}

\newcommand{\M}{{\mathcal M}}
\newcommand{\ra}{\rightarrow}
\newcommand{\expec}{{\mathbb E \,}}
\newcommand{\torus}{\T}
\newcommand{\Mjutska}{K}
\newcommand{\beqno}{\begin{eqnarray*}}
\newcommand{\eeqno}{\end{eqnarray*}}

\newcommand{\real}{\R}
\newcommand{\refeq}[1]{(\ref{#1})}

\def\ds{\displaystyle}

\newcommand{\beqla}[1] {\begin {eqnarray}\label{#1}}

\def \e {\varepsilon}
\def \UCM {\mathbf{u}}
\def \beq {\begin {eqnarray}}
\def \eeq {\end {eqnarray}}
\def \ba {\begin {eqnarray*}}
\def \ea {\end  {eqnarray*}}

\def \bra{\langle}
\def \cet{\rangle}
\def \prob{{\mathbb P}}

\def \Noise {\mathcal{E}}
\def \noise {\varepsilon}
\def \prior {\Pi}

\begin{document}

\begin{center}
{\LARGE Discretization-invariant Bayesian inversion

\vspace{.3cm}

and Besov space priors}

\vspace{1cm}

{\sc Matti Lassas}\\
Department of Mathematics and Statistics\\
University of Helsinki, Helsinki, Finland\\

\vspace{.6cm}

{\sc Eero Saksman}\\
Department of Mathematics and Statistics\\
University of Helsinki, Helsinki, Finland\\

\vspace{.6cm}

{\sc  Samuli Siltanen}\\
Tampere University of Technology, Department of Mathematics \\
P.O.Box 553, 33101 Tampere, Finland

\end{center}

\noindent {{\sc Abstract.}\it
 {\mntext Bayesian solution of an inverse problem for} indirect measurement $M = AU + \Noise$ is considered,
{\mntext  where $U$ is a function on a domain of $\R^d$.}
Here $A$ is a smoothing linear operator and $\Noise$ is {\ntext Gaussian} white noise.
{\mntext The data is} a realization $m_k$ of the
random variable $M_k = P_kA U+P_k\Noise$, where $P_k$ is a linear, finite dimensional operator related to measurement device.
To allow computerized inversion,
the unknown is discretized as $U_n=T_nU$, {\mntext where $T_n$ is a finite dimensional projection,}
leading to the computational measurement model $M_{kn}=P_k A U_n + P_k\Noise$.
Bayes formula gives then the posterior distribution
$\pi_{kn}(u_n\,|\,m_{kn})\sim\prior_n(u_n)\,\exp(-\frac{1}{2}\|m_{kn} - P_kA u_n\|_2^2)$ in $\R^d$,
and the mean $\UCM_{kn}:=\int
u_n \, \pi_{kn}(u_n\,|\,m_k)\, du_n$ is considered as the reconstruction of $U$. We discuss a systematic way of choosing
prior distributions $\prior_n$ for all $n\geq n_0>0$ by achieving them as projections of a distribution in a infinite-dimensional limit case.
Such choice of prior distributions is {\em discretization-invariant} in the sense that $\prior_n$ represent the same {\em a priori} information for all $n$ and that the mean $\UCM_{kn}$ converges to a limit estimate as $k,n\ra\infty$.
Gaussian smoothness priors and wavelet-based Besov space priors are shown to be discretization invariant. In particular, Bayesian inversion in dimension two with $B^1_{11}$ prior is related to penalizing the $\ell^1$ norm of the wavelet coefficients of $U$.}

 {\bf Keywords.} Inverse problem, statistical inversion, Bayesian inversion, discretization invariance,
reconstruction, wavelet, Besov space

 {\bf AMS classification.} 60F17, 65C20, 42C40


\vspace{1cm}

\section{Introduction}

Consider a quantity $U$ that can be observed indirectly through a relation
\begin{eqnarray} \label{basicequation A1}
  M = AU + \Noise,
\end{eqnarray}
where $A$ is a smoothing linear operator and $\Noise$ is white noise. Here $U$ and $M$ are considered as continuum objects, or functions defined on subsets of $\R^d$, so that our discussion applies to classical models of mathematical physics such as Laplace, Maxwell, Helmholtz or Schr\"odinger equation. We are interested in the use of Bayesian inversion to find information about $U$ from measurement data concerning $M$. Let  $U(x,\omega),\ M(y,\omega)$ and $\Noise(y,\omega)$ be random functions where $\omega\in \Omega$ is an element of a complete probability space  $(\Omega,\Sigma,\prob)$ and $x$ and $y$ denote the variables in domains of Euclidean spaces. We analyze Bayesian estimates of $U$ when a continuum model of the form (\ref{basicequation A1}) is approximated by finite-dimensional models to allow computerized inversion.

Assume that a measurement device provides us with a realization of the random variable
\begin{equation}\label{eq: measmodel_k k}
  M_k = P_kM = A_k U+\Noise_k,
\end{equation}
where $A_k=P_k A$ and $\Noise_k=P_k\Noise$. Here $P_k$ is a linear operator related to the device; for simplicity we assume that $P_k$ is an orthogonal projection with $k$-dimensional range. We call (\ref{eq: measmodel_k k}) the {\it practical measurement model}
and (\ref{basicequation A1}) the {\it continuum model}. Realizations of measurements are denoted by $m_k=M_k(\omega_0)$ and $m=M(\omega_0)$,
respectively, where $\omega_0\in \Omega$ is a specific element in the probability space. 

This study concentrates on the inverse problem
\begin{eqnarray}\label{inverseproblem}
  \mbox{given a realization } M_k(\omega_0)\mbox{, estimate }U,
\end{eqnarray}
where the estimates in question are means and confidence intervals related to a Bayesian posterior probability distribution.

For example, consider the brain imaging method called magnetoencephalography ({\sc meg}),
see e.g. \cite{HariLounasmaa}. Maxwell's equations describe how synchronized neuronal currents $U$ in the cerebral cortex produce a magnetic field $AU$ that can be measured at the surface of the head. Let $\Noise$ denote the magnetic fields produced by all external sources; then the continuum model (\ref{basicequation A1}) describes the total magnetic field $M=AU+\Noise$. In practice one measures the inner products $\bra M,\phi_j\cet$, $j=1,2,\dots,k,$ where $\phi_j$ are linearly independent device functions corresponding to measuring the flux of the magnetic field $M$ through a small surface determined by the $j$th measurement unit ({\sc squid}).
As an idealization, let us assume that $\phi_j$ are orthonormal so that $P_kv=\sum_{j=1}^k \bra v,\phi_j\cet\phi_j$ is an orthogonal projection.
Then {\sc meg} data is modelled by $P_kM$.

Computational solution of (\ref{inverseproblem}) using Bayesian inversion involves discretization of the unknown quantity $U$. We assume that $U$ is {\em a priori} known to take values in a function space $Y$. Choose a linear projection $T_n:Y\ra Y$ with $n$-dimensional range $Y_n$, and define a random variable $U_n:=T_n U$ taking values in $Y_n$. This leads to the {\it computational model}
\begin{equation}\label{eq: measmodel_eka_kn}
  M_{kn} = A_k U_n+\Noise_k
\end{equation}
involving two {\em independent} discretizations: $P_k$ is related to the measurement device and $T_n$ to finite representation of the unknown.

In the above application related to {\sc meg}, the projection $T_n$ corresponds to an approximate representation of the electromagnetic sources in the brain using a finite set of basis functions (defined for instance according to a finite element method).

Note that the model (\ref{eq: measmodel_eka_kn}) is virtual in the sense that $U_n$ appears neither in the continuum model (\ref{basicequation A1}) nor in the practical measurement model (\ref{eq: measmodel_k k}). In particular, measurement  ${M}_{kn}(\omega_0)$ is related to the computational model but not to the practical measurement model. This is why we use $m_k=M_k(\omega_0)$ as the given data.

Denote the probability density function of the random variable $M_{kn}$ by $\Upsilon_{kn}(m_{kn})$. The {\em posterior density} for $U_n$ is given by the Bayes formula:
\begin{eqnarray}\label{posterior2}
  \pi_{kn}(u_n\,|\,m_{kn})
  =
  \frac{\prior_n(u_n)\,\exp(-\frac{1}{2}\|m_{kn} - A_k u_n\|_2^2)}{\Upsilon_{kn}(m_{kn})},
\end{eqnarray}
where the exponential function corresponds to (\ref{eq: measmodel_eka_kn}) with white noise statistics
{\ntext with identity variance}, and {\em a priori} information about $U$ is expressed in the form of
a {\em prior density} $\prior_n$ for the random variable $U_n$. The density $\prior_n$ assigns high probability to functions that are
typical in light of {\em a priori} information, and low probability to atypical functions.

We can now state the inverse problem (\ref{inverseproblem}) more specifically:
\begin{eqnarray}\label{inverseproblem3}
  \mbox{given a realization }m_k= M_k(\omega_0)\mbox{, estimate }U\mbox{ by }\UCM_{kn},
\end{eqnarray}
where the {\ntext conditional mean ({\sc cm})} estimate (or posterior mean estimate) $\UCM_{kn}$ is
\begin{equation}\label{meankaava2}
  \UCM_{kn}:=\int_{Y_n} u_n \, \pi_{kn}(u_n\,|\,m_k)\, du_n.
\end{equation}
Note that formula (\ref{meankaava2}) differs from the conventional definition of posterior mean estimate since it involves $m_k= M_k(\omega_0)$
 instead of $m_{kn}= M_{kn}(\omega_0)$.

The estimate $\UCM_{kn}$ and confidence intervals for it are typically computed approximately with Markov chain Monte Carlo ({\sc mcmc}) methods involving simulation software for the finite model (\ref{eq: measmodel_eka_kn}). The solution strategy (\ref{inverseproblem3}) has been applied to image restoration \cite{besag}, geological prospecting \cite{Nicholls,andersen}, atmospheric and ionospheric remote sensing \cite{barbara, tamminen, nygren, markkanen, lehtinen2}, medical X-ray tomography \cite{dentalpaper1,dentalpaper2,hanson1} and electrical impedance imaging \cite{Kaipio1,andersen2}. For general reference on Bayesian inversion see \cite{KS,mosegaard,evansstark,tarantolabook,EkiDanielaBook}.

We remark that applying the {\sc mcmc} solution strategy requires also practical implementation of the operator $A_k$. In the case of {\sc meg} imaging $A_k$ corresponds to computing the electromagnetic field $A_kU_n$ using the discrete current $U_n$ and an approximate numerical solution to Maxwell's equations. We neglect the effects of numerical error in the implementation of $A_k$ in this paper. (One possibility to take this error into account is to use the approximation error model of Kaipio and Somersalo \cite[Section 5.8]{KS}, but this leads to non-Gaussian noise statistics in general and thus falls outside the scope of this discussion.)

Summarizing, our starting point is the infinite-dimensional continuum model (\ref{basicequation A1}). A measurement instrument provides us with finite-dimensional noisy data $m_k=M_k(\omega_0)$ described by the practical measurement model (\ref{eq: measmodel_k k}). Our aim is to use $m_k$ to find information about the unknown $U$. To allow computerized inversion we construct the fully discrete computational model (\ref{eq: measmodel_eka_kn}) involving {\em a priori} information about $U$, and we write down the posterior distribution (\ref{posterior2}). Finally, we use formula (\ref{meankaava2}) and numerical methods to estimate $U$ with $\UCM_{kn}$. However,  in definition (\ref{meankaava2}) the data $m_k$ comes from the practical measurement model (\ref{eq: measmodel_k k}), while $\pi_{kn}$ is related to the computational model (\ref{eq: measmodel_eka_kn}). Taking this incompatibility into account is one of the central novelties in this work.

Constructing $T_n$ and $\prior_n$ is the core difficulty in Bayesian inversion. Often there is no natural discretization for the continuum quantity $U$, so $n$ can be freely chosen. Consequently, $T_n$ and $\prior_n$ should in principle be described for all $n>0$, or at least for an infinite sequence of increasing values of $n$. Also, updating our measurement device may increase $k$ independently of $n$. This work is motivated by the need to avoid the following unwanted phenomena:
\begin{itemize}
\item[(a)] {\bf The estimates $\UCM_{kn}$ diverge as $n\ra\infty$.} In this case investing more computational resources to modeling our unknown does not necessarily result in improved reconstructions.
\item[(b)] {\bf The estimates $\UCM_{kn}$ diverge as $k\ra\infty$.} In this case performing more measurements may lead to worse reconstructions.
\item[(c)] {\bf Representation of {\em a priori} knowledge is incompatible with discretization.} It is reported in \cite{LS} that discrete (non-Gaussian) total variation priors converge to a Gaussian smoothness prior as $n\ra\infty$. In this case one makes the mistake of specifying different {\em a priori} information for different values of $n$. See Appendix \ref{sec:discretizationexamples} below for more details.
\end{itemize}
A choice of $T_n$ and $\prior_n$ is called {\em discretization-invariant} if it avoids (a)--(c).

We construct prior distributions for $U$ in the infinite-dimensional space $Y$. Then the random variable $U_n=T_nU$ takes values in the finite-dimensional subspace $Y_n\subset Y$ and represents approximately the same {\em a priori} knowledge {\ntext as} $U$. Further, we analyze convergence of $\UCM_{kn}$ using a deterministic function called {\em reconstructor} that almost surely maps a given measurement to the conditional mean estimate. For example, the reconstructor $\mathcal{R}_{M_{kn}}(U_n|\,\,\cdotp\,)$ corresponding to the computational model (\ref{eq: measmodel_eka_kn}) takes the measurement data $m_k=M_k(\omega_0)$ to the mean $\UCM_{kn}$ defined in (\ref{meankaava2}). The infinite-dimensional model $M = AU + \Noise$ has a reconstructor $\mathcal{R}_M(U|\,\,\cdotp\,)$ as well. Theorem \ref{Corollary-teoreema} below states under suitable assumptions on $U,T_n$ and $P_k$ that
\begin{equation}\label{introrecconv}
 \lim_{n,k\to \infty}\mathcal{R}_{M_{kn}}(U_n | m_k ) =\mathcal{R}_M(U|m),
\end{equation}
where the realization $m_k=P_km$  comes from a realization $m=M(\omega_0)$ of the random variable $M$.

Our proving strategy involves another measurement model analogous to (\ref{eq: measmodel_eka_kn}):
\begin{equation}\label{eq: measmodel_kn}
  \Meas_{kn} = A_k U_n+\Noise,
\end{equation}
where the noise is not finite-dimensional. The {\ntext noise} models (\ref{basicequation A1}) and (\ref{eq: measmodel_kn}) now contain the same (continuum) white noise process. This allows us to prove convergence results in a fixed function space. Theorem  \ref{convergence -theorem} below states under suitable assumptions on $U,T_n$ and $P_k$ that if $\lim_{k\ra \infty}m_k=m$ then
$$
 \lim_{n,k\to \infty}\mathcal{R}_{\Meas_{kn}}(U_n | m_k ) =\mathcal{R}_M(U|m).
$$
Formula (\ref{introrecconv}) follows by showing that the reconstructors coincide: $\mathcal{R}_{\Theta_{kn}}(U_n|m_k) = \mathcal{R}_{M_{kn}}(U_n|m_k)$ for $m_k\in \hbox{Ran}\,(P_k)$.
{\nnntext We will consider also more general reconstructors $\mathcal{R}_{M}(g(U)|m)$ that can be used
to analyze convergence of confidence intervals.}

Our way of using infinite-dimensional limit processes is one possibility to achieve discretization-invariance since we avoid problems (a)--(c).

We are especially interested in discretization-invariant and edge-preserving Bayesian inversion. Total variation regularization of Rudin, Osher and Fatemi \cite{OsherRudin} penalizes the $L^1$ norm of derivatives and yields edge-preserving reconstructions in practice
\cite{dobson,vogel,Kaipio1,dentalpaper1,dentalpaper2}. These results are equivalent to computing {\em maximum a posteriori}  estimates {\ntext using a total variation prior and a Gaussian likelihood}.  However, Bayesian inversion with {\mntext discretized} total variation prior distribution,
\begin{eqnarray}\label{TV prior}
  \prior_{U_n}(u)= c_n\exp(-\alpha_n\|\nabla u\|_{L^1(0,1)}),\quad u\in Y_n,
\end{eqnarray}
{\mntext where $Y_n\subset W^{1,1}_0(0,1)$ are finite dimensional subspaces},
is not discretization-invariant \cite{LS}; in particular, conditional mean estimates
({\ntext i.e.\ posterior mean estimates)} lose their edge-preserving quality as $n\ra\infty$.

Wavelet-based Bayesian inversion using Besov space priors is used in applications with results similar to total variation regularization,
see \cite{BayesWavelet,maaria,kati,multipatent}. The Besov space $B_{11}^1$ that bounds $L^1$ norms of (band-filtered) first derivatives
similarly to (\ref{TV prior}) is useful for image processing, see Meyer \cite{meyer}. One of the main results of this paper is that Besov
priors are discretization-invariant; we emphasize that this gives us non-Gaussian discretization-invariant priors.
The proof is based on the analysis of the infinite-dimensional limit case and includes a quantitative estimate on the speed of convergence
of reconstructors. Further, we show in Section \ref{sec:Besov_deblur} for the special case of two-dimensional deblurring that $B_{11}^1$ inversion using $\UCM_{kn}$ reduces to applying {\ntext $\ell^1$-type prior} on the wavelet coefficients of $U_n$---a combination of two well-known computational methods. See Section \ref{sec:Besov_deblur} for more details. There is an
interesting parallel to algorithms computing {\sc map} estimates with penalty on the $\ell^1$ norm of Fourier or wavelet coefficients \cite{donoho1,donoho2,Daubechies2,loris}.

Let us review previous literature on the topic. The study of Bayesian inversion in infinite-dimensional function spaces was initiated by
Franklin \cite{franklin} and continued by Mandelbaum \cite{mandelbaum}, Lehtinen, P\"aiv\"arinta and Somersalo \cite{LPS}, Fitzpatrick
\cite{fitzpatrick}, and Luschgy \cite{luschgy}.  The concept of discretization invariance was formulated by Markku Lehtinen in the 1990's
and has been studied by Lasanen \cite{La}, Piiroinen \cite{Pi}, D'Ambrogi, M\"aenp\"a\"a and Markkanen \cite{barbara}, and Lasanen and
Roininen \cite{SariLassi}. A definition of discretization invariance similar to the above was given in \cite{LS}. For other kinds of
discretization of continuum objects in the Bayesian framework, see \cite{hanson1,Nicholls}. The use of wavelets and Besov spaces in
statistical algorithms is discussed in \cite{ASS,AS,BayesWavelet,leporini,choi,waveletstatbook}. For regularization based approaches for
statistical inverse problems, see \cite{engl1,engl2,hohage,Pikkarainen}. The relationship between continuous and discrete (non-statistical)
inversion is studied in Hilbert spaces in \cite{Vogel2}. See \cite{Borcea} for specialized discretizations for inverse problems.

We remark that working entirely within the computational model (\ref{eq: measmodel_eka_kn}) would require using a realization $M_{kn}(\omega_0)$ of the random variable $M_{kn}$ instead
of $m_k=M_k(\omega_0)$ in formula (\ref{meankaava2}). The starting point of inversion in earlier studies of discretization invariance is indeed
the random variable $M_{kn}$ or its realization. 
However, the appropriate model of data 
given by the measurement device is a realization 
$m_k$ related to model (\ref{eq: measmodel_k k}).
 Rigorous analysis of this incompatibility 
is a central novelty in this paper. 
To emphasize this aspect we show that inversion from $m_k$ using Gaussian smoothness priors or Besov space priors is discretization-invariant. {\ntext We do not discuss non-Gaussian noise models in this paper.}

This paper is organized as follows. In Section \ref{sec:deblurring} we discuss discretization-invariant Bayesian inversion using Gaussian and Besov priors. In Section \ref{sec:general} we present the general theory of discretization invariance.
In Section \ref{sec:randomBesov} we discuss random variables with values in Besov spaces. After proving some technical results in Section
\ref{sec:technical} we apply them in Section \ref{sec:BesovConv} to prove convergence of reconstructors arising from Besov priors. {\nnntext In Appendix \ref{sec:discretizationexamples} we consider examples.}

Below $\langle\cdot,\cdot\rangle$ refers to pairing of either generalized functions with test functions, or a Banach space with its dual.
We denote by $\langle\cdot,\cdot\rangle_X$ the inner product in a Hilbert space $X$. The notation $L(X,Y)$ stands for the space of bounded
linear operators between Banach spaces $X$ and $Y$, and $L(X,X)$ is abbreviated as $L(X)$. We occasionally denote the norm $\|\,\cdot\,\|_{L(X,Y)}$ by $\|\,\cdot\,\|_{X\ra Y}$. The specific element $\omega_0\in\Omega$ of the probability space denotes realizations of measurements throughout the paper.

{\bf Acknowledgements:} the authors thank Markku Lehtinen, Heikki Haario and Marko Laine for useful discussions and the anonymous referees 
whose thoughtful 
comments helped us to improve the paper. ML and SS were supported
by National Technology Agency of Finland (TEKES) under Contract 206/03,
Finnish Centre of Excellence in Inverse Problems Research and
PaloDEx Group.
ES was supported by the Academy of
Finland, projects no. 113826 and 118765, and by the Finnish
Center of Excellence in Analysis and Dynamics.

\section{Example: discretization-invariant deblurring}\label{sec:deblurring}

{\mntext In this section we give examples of discretization-invariant prior distributions and consider
a simple inverse problem to give a flavour of our results to the reader. The precise definitions
in a more general setting are postponed to the later sections.}

We discuss Bayesian deconvolution using a Gaussian smoothness prior and a Besov space prior. In both cases we define the prior distribution
in the continuous context and then marginalize it to discrete cases to allow practical computation. The Gaussian case is shown to be
related to deterministic Tikhonov regularization with derivative penalty.

For simplicity, we ignore boundary effects by considering periodic functions. The loss of generality is not too bad from the practical
point of view since the periodic analysis covers compactly supported non-periodic cases.

Let $\T^2$ be the two-dimensional torus constructed by identifying parallel sides of the square $D=(0,1)^2\subset \R^2$; we model periodic
images as elements of function spaces over $\torus^2$. The {\mntext continuum model} is $M=AU+\Noise$ with convolution operator $A$
defined by
\begin{equation}\label{A_conv}
  Au(x)=\int_{\T^2} \Phi(x-y)u(y)\,dy,
\end{equation}
where $\Phi\in C^\infty(\T^2)$ is a point spread function.

\subsection{Gaussian smoothness prior}\label{sec:Gaussianexample}

For any  $s\in\R$, let $H^s(\T^2)$ be the $L^2$-based Sobolev space equipped with Hilbert space inner product
\begin{eqnarray}\label{eq: inner product in Hs}
\bra \phi,\psi\cet_{H^s(\T^2)}=\int_{\T^2} ((I-\Delta)^{s/2}\phi)(x)\, ((I-\Delta)^{s/2}\psi)(x)\,dx.
\end{eqnarray}
Note that $H^0(\T^2)=L^2(\T^2)$.

Recall that a generalized Gaussian random variable $V$ takes values in the space of generalized functions, and the pairing $\bra
V,\phi\cet$ with any test function $\phi\in C^\infty(\T^2)$ is a Gaussian random variable taking values in $\R$, see \cite{roz}. {\mntext The Gaussian random variables we will consider below are assumed to
take values in some Hilbert space, typically in a Sobolev space $H^s(\T^2)$, where the smoothness
index $s\in \R$ may also be negative.
Now,} if $V$
takes values in a Hilbert space $X$ we say that $V$ has the covariance operator $C_V:X\to X$ if
\begin{equation}\label{covop}
  \expec (\bra V-\expec V,\phi\cet_X\, \bra V-\expec V,\psi \cet_X)= \bra C_V\phi,\psi\cet_{X},
\end{equation}
with any $\phi,\psi\in X$. Here $\langle\cdot,\cdot\rangle_X$ stands for the inner product in $X$.

{\mntext Next we analyze a simple measurement model as an example.}
Consider the {\mntext continuum} model $M=AU+\Noise$ where the convolution operator (\ref{A_conv}) is now viewed as a smoothing map
$A:H^s(\T^2)\to C^\infty(\T^2)$.

We construct the smoothness prior by choosing $U$ to be a generalized Gaussian random variable taking values in $H^{-1}(\T^2)$ and having
expectation $\expec U=0$ and covariance operator $C_U=\alpha^{-1}(I-\Delta)^{-2}$. Here $\alpha>0$ is a regularization parameter. The
operator $C_U$ corresponds formally to the prior
\begin{eqnarray*}
 \pi_U(u)\underset{\mbox{\tiny \em formally}}=  c\exp(- \frac{\alpha}{2} \|u\|_{H^1(\torus^2)}^2)
\end{eqnarray*}
and generates the discrete smoothness priors widely used in practice. However, it is curious that in spite of the term {\em smoothness
prior} the realizations of $U$ are almost surely not even $L^2(\T^2)$ functions, let alone differentiable. This is why we need to consider
$U$ as taking values in some space $H^s(\T^2)$ with negative smoothness index $s$; the value $s=-1$ is chosen just for convenience.

White noise $\Noise$ is a generalized Gaussian random variable with expectation $\expec \Noise=0$. Let us discuss the choice of an
appropriate covariance operator: {\mntext the standard definition of white noise as a generalized
random variable is that}
\beq\label{eq: covariance for white noise}
  \expec(\langle\Noise,\phi\rangle\langle\Noise,\psi\rangle)=\langle\phi,\psi\rangle,\qquad \mbox{for all }\phi,\psi\in C^\infty(\torus^2),
\eeq
{\mntext where $\langle\cdot,\cdot\rangle$ denotes the pairing between generalized functions and test functions.
To consider the white noise $\Noise$ as a Hilbert-space-valued random variable, we
can choose the Hilbert space to be any Sobolev space $H^s(\torus^2)$, $s<-1$.
 One {\mntext possible} choice is
to consider $\Noise$ as taking values in $H^{-2}(\T^2)$ and choose the covariance operator $C_\Noise=(I-\Delta)^{-2}:H^{-2}(\T^2)\ra
H^{-2}(\T^2)$ as defined in (\ref{covop}). {\mntext Note that
realizations of $\Noise$ belong to  $L^2(\torus^2)$ only with probability
zero, and this is why we need to use Sobolev spaces with negative
smoothness indices $s$.}}

Now the continuous framework for inversion is in place. Let $m=M(\omega_0)$ be a realization of the measurement $M=AU+\Noise$. Since both the
prior and noise statistics are Gaussian, the conditional mean estimate coincides with the location of the maximum of the posterior density.
Thus we can evaluate the {\sc cm} estimate $\UCM$ as
\begin{equation}\label{ekaGCM}
  \UCM = \mbox{arg}\!\!\!\!\!\max_{u\in H^1(\T^2)}
  \big\{\exp(-\frac 12 \|A u\|_{L^2(\T^2)}^2 + \langle m,A u \rangle -\langle C_U^{-1} u,\,u\rangle_{H^{-1}(\T^2)}) \big\}.
\end{equation}
We omitted in formula (\ref{ekaGCM}) the constant term $\|m\|_{L^2(\T^2)}^2$ in the formal expansion
\begin{equation}\label{apukaava}
  \|m-A u\|_{L^2(\T^2)}^2 = \|A u\|_{L^2(\T^2)}^2 -2 \langle m, A u \rangle_{L^2(\T^2)} +\|m\|_{L^2(\T^2)}^2,
\end{equation}
which is well-defined only when $m\in L^2(\T^2)$, and this happens with probability zero. Also, the smoothing properties of the
operator $A$ make it possible to replace the formal quantity $\langle m, A u \rangle_{L^2(\T^2)}$ in (\ref{apukaava}) by the
rigorously defined pairing $\langle m,A u \rangle$ in (\ref{ekaGCM}). Note that we can write (\ref{ekaGCM}) in the form
$$
  \UCM = \mbox{arg}\!\!\!\!\!\min_{u\in H^1(\T^2)}
  \big\{\frac 12 \|A u\|_{L^2(\T^2)}^2 - \langle m,A u \rangle +\frac{\alpha}{2}\|(I-\Delta)^{1/2}u\|_{L^2(\T^2)}^2 \big\}.
$$

{\ntext The} practical measurement model is $M_k = P_kAU+P_k\Noise$, where $P_k:L^2(\T^2)\to L^2(\T^2)$ is an
orthogonal projection with $k$-dimensional range. We require that the sequence $P_k$ converges 
strongly to the identity operator on $L^2(\T^2)$ as $k\to \infty$. For example, $P_k$ may measure averages of $AU$ over $k$ pixels on $\T^2$ and
construct a piecewise function or compute a truncated Fourier series expansion to get an element of $L^2(\T^2)$.  Here $k$ can be arbitrarily large, enabling models of imaging devices
with any resolution.

Let us now turn to practical inversion {\mntext where all appearing quantities are finite dimensional}.
We need to discretize the unknown. Take $T_n:H^{-1}(\T^2)\to H^{-1}(\T^2)$ to be truncations of Fourier series expansions to $n$ lowest frequency terms. Then $T_n$ are linear orthogonal projections with $n$-dimensional range $Y_n$ converging strongly to the identity operator on $H^{-1}(\T^2)$ as $n\to \infty$. Let $U$ be as in the {\mntext continuum} case above and define a random variable $U_n:=T_n U$ taking values in $Y_n$. The conditional mean estimate for $U_n$ (determined using the posterior distribution corresponding to model $M_{kn}=P_kAU_n + \Noise_k$) is
\begin{eqnarray}
  \UCM_{kn}
  &:=& \label{Gauss_deblur_discrete}
  \mbox{arg}\!\!\!\!\!\!\!\!\min_{u\in Y_n\cap H^1(\T^2)}
  \big\{\frac 12 \|P_k A T_n u\|_{L^2(\T^2)}^2 - \langle m_k,P_k A T_n u \rangle\\
  &&\nonumber
  \qquad \qquad\qquad +\frac{\alpha}{2}\langle (T_n C_U T_n)^{-1} u, \,u\rangle_{H^{-1}(\T^2)} \big\},
\end{eqnarray}
where $(T_n C_U T_n)^{-1}$ is the inverse of $T_n C_U T_n:Y_n\to Y_n$.
{\nnntext Here, $m_k=M_k(\omega_0)$ is the realization of measurement
$M_k$ in the practical measurement model.}

Theorem \ref{convergence -theorem} below implies
in particular that $\UCM_{kn}\ra \UCM$  as $k,n\ra\infty$, showing that Bayesian deblurring with Gaussian smoothness prior is discretization-invariant.

Much of the material in this Subsection \ref{sec:Gaussianexample} is due to Lasanen \cite{La} and Piiroinen \cite{Pi}. We emphasize that we assume given a realization $m_k=M_k(\omega_0)$ from practical measurement model (\ref{eq: measmodel_k k}) and that we do not have available any realization of $M_{kn}$.  Because of this, $m_k$ is used in formula (\ref{Gauss_deblur_discrete}). This is the main novelty of our Gaussian example compared to \cite{La}.

We close this section by discussing a connection to generalized Tikhonov regularization, here defined as finding the element in $u\in Y_n\cap
H^1(\T^2)$ that minimizes the functional
$$
  \frac 12 \|m_k-P_kA u\|_{L^2(\T^2)}^2+\frac{\alpha}{2}\|u\|_{L^2(\T^2)}^2+\frac{\alpha}{2}\|\nabla u\|_{L^2(\T^2)}^2.
$$
After similar modification as above we can define the Tikhonov solution by
\begin{eqnarray}
  \mathbf{u}_{kn}^T
  &=& \label{tikhonov}
  \mbox{arg}\!\!\!\!\!\!\!\!\min_{u\in Y_n\cap H^1(\T^2)} \Big\{\frac 12 \|P_kA u\|_{L^2(\T^2)}^2 - \langle m_k,P_k A u \rangle\\
  && \nonumber
  \qquad \qquad\qquad +\frac{\alpha}{2}\langle (I-\Delta)u,\, u\rangle_{L^2(\T^2)} \Big\}.
\end{eqnarray}
The quadratic forms defined in (\ref{Gauss_deblur_discrete}) and (\ref{tikhonov}) are the same, so $\UCM_{kn}=\UCM_{kn}^T$. Thus our
results apply to the convergence analysis of Tikhonov regularization as well.

\subsection{Besov prior}\label{sec:Besov_deblur}

Gaussian smoothness priors are designed for representing the prior information that the unknown physical quantity $U$ does not vary
sharply. However, sometimes we know {\em a priori} that $U$ is piecewise regular, and other kind of priors are needed. One could use
Gaussian hyperpriors as in \cite{Kaipio2,EkiDaniela2} or geometric priors as in \cite{Nicholls}. We discuss a different approach that
allows analysis of an infinite-dimensional limit case and consequently leads to discretization-invariant inversion.

We replace the discretization-dependent \cite {LS} total variation prior
$$
  \prior(u)\underset{\mbox{\tiny \em formally}}= c\exp(-\alpha\|\nabla u\|_{L^1})
$$
by the discretization-invariant Besov space prior
$$
  \prior(u)\underset{\mbox{\tiny \em formally}}= c\exp(-\alpha\| u\|_{B_{11}^1(\torus^2)}).
$$
Since the norm of $B_{11}^1(\torus^2)$ bounds the $L^1$ norms of (band-limited) first derivatives of $u$, we expect that the $B_{11}^1$
prior can be used for edge-preserving inversion. It is computationally convenient that $B_{11}^1$ functions can be written {\ntext using a compactly supported} wavelet bases
and the $B_{11}^1$ norm can be computed as a weighted sum of wavelet coefficients, see  Appendix \ref{app:wavelet}.

The {\mntext continuum} measurement model is $M=AU+\Noise$ where the linear convolution operator (\ref{A_conv}) is considered as a smoothing map
$A:B_{11}^1(\torus^2)\ra C^\infty(\torus^2)$.

We construct the $B_{11}^1(\torus^2)$ prior by expanding $U$ in the wavelet basis as the sum
\begin{eqnarray*}
  U = \sum_{\ell= 1}^\infty X_{\ell}\psi_{\ell}
\end{eqnarray*}
with each random coefficient $X_\ell$ distributed independently {\ntext according to $\pi_X(x)=c\exp(-|x|)$,
where $c=2e^{-1}$ is the normalization constant.} This distribution arises naturally due to
the wavelet characterization of Besov spaces. The scale of the wavelet basis
functions becomes finer when $\ell$ {\ntext increases}; for the exact bookkeeping of scales and locations of $\psi_{\ell}$ as function of $\ell$ we
refer to Appendix \ref{app:wavelet}.

We take $T_n:B_{11}^1(\torus^2)\to B_{11}^1(\torus^2)$ to be the finite-dimensional projections
\begin{eqnarray}\label{eq: cut of series}
  T_n\left(\sum_{\ell=1}^\infty c_\ell \psi_\ell\right) =  \sum_{\ell=1}^n c_\ell \psi_\ell
\end{eqnarray}
that simply truncate the wavelet expansion to $n$ first terms. These operators $T_n$ converge strongly to
{\ntext the identity as} $n\to
\infty$. {\ntext For each $n$,} define a random variable $U_n:=T_n U$ taking
values in $Y_n:=\mbox{span}(\psi_1,\dots,\psi_n)$ and consider the model $M_{kn}=P_kAT_nU + \Noise_k$
{\mntext where $P_k:L^2(\T^2)\to L^2(\T^2)$ is an
orthogonal projection with $k$-dimensional range.}

With the above choices the posterior distribution of $U_n$ takes the computationally efficient form in terms of the {\mntext wavelet coefficients $X^{n}:=(X_1,X_2,\dots,X_n)$ of $U_n$, namely, the probability 
distribution of $X^n$ conditioned with $M_{kn}$ is} 
\begin{equation}\label{Besovposterior}
\pi_{X^n|M_{kn}}(x_1,\dots, x_n,m_k)=
  C\exp\Big( - \frac 12 \big\|m_k - A\sum_{\ell=1}^n x_\ell\psi_\ell\big\|_{L^2(\torus^2)}^2 -\alpha\sum_{\ell=1}^n |x_\ell|\Big).
\end{equation}
The conditional mean estimate can be computed approximately e.g.\ using a Markov chain Monte Carlo algorithm such as the Gibbs sampler or the
Metropolis-Hastings method. Sampling from $L^1$ distributions is explained e.g. in
\cite{nenonen,Kaipio1} and \cite[Section 3.3.2]{KS}, and modifying those methods for computing the mean of (\ref{Besovposterior}) involves only adding the fast wavelet transform appropriately.

Our new results below show that Besov priors are discretization-invariant. First, by Theorems \ref{Corollary-teoreema} and \ref{th:besovconvergence} and Corollary \ref{co:10} the estimates $\UCM_{kn}$ converge when $n,k\ra\infty$ either separately or simultaneously. Thus Besov priors do not involve the possible difficulties (a) and (b) mentioned in the introduction. Second, the prior distributions defined in $Y_n$ converge to a limit distribution in $Y$ as $n\ra\infty$ (Proposition \ref{discret inv}). Thus we do not have the unwanted phenomenon (c) of the introduction.

\section{General theory of discretization invariance}\label{sec:general}

Consider two independent random variables $U$ and $\Noise$ and the measurement model $M=AU+\Noise$. We construct a rigorous stochastic
framework for discretization-invariant Bayesian inversion using the following diagram of spaces and maps:
\begin{equation}\label{basicdiagram}
\begin{array}{cccccc}
&&&&&\\
 Y &  \stackrel{A}{\longrightarrow}
& \ S^{1}& \ \subset\quad S^{1/2}\!=\!Z\quad \subset\quad S & \\
\rotatebox{90}{$\in$} &&&&\!\!\!\!\!\!\!\!\!\!\!\!\!\!\!\rotatebox{90}{$\in$}\\
U(\omega_1) &&&& \!\!\!\!\!\!\!\!\!\!\!\!\!\!\!\Noise (\omega_2)
\\ &&&&
\end{array}
\end{equation}
\smallskip
\noindent Our immediate task is to define  all the objects in (\ref{basicdiagram}).

Let $(\Omega,\Sigma,\prob)$ a complete probability space with product structure: $\Omega=\Omega_1\times\Omega_2$ and $\Sigma=\overline{\Sigma_1\otimes \Sigma_2}$ and $\prob=\prob_1\otimes \prob_2$.

Recall the notion of a random variable taking values in a Banach space. {\ntext For {\nnntext any given} separable Banach space $X$
we denote} the dual of $X$ by $X'$. Let
$\B_X$ be the Borel $\sigma$-algebra of $(X,\tau^w)$ with $\tau^w$ the weak topology of $X$. Note that the separability of $X$ verifies that $\B_X$ coincides with the Borel
$\sigma$-algebra of $(X,\tau^n)$, where $\tau^n$ is the norm topology of $X$.  An $X$-valued
random variable $V$ is simply any measurable map $V:(\Omega,\Sigma)\to (X,\B_X)$. In this paper we consider
only random variables with values in separable Banach spaces.

We now assume that the space  $Y$ in the diagram (\ref{basicdiagram}) is a separable Banach space and that $U=U(\omega_1)$ is a $Y$-valued random variable.

The measurement noise $\Noise=\Noise(\omega_2)$ with $\omega_2\in\Omega_2$ is a Gaussian random variable taking values in a
separable Hilbert space $S$; the expectation satisfies $\expec \Noise=0$, and the covariance operator $C_\Noise:S\ra S$, is defined by requiring
$$
\expec
(\bra s_1,\Noise\cet_S\bra \Noise,s_2\cet_S)=\bra s_1,C_\Noise s_2\cet_S
\quad \hbox{ for }s_1,s_2\in S.
$$
We assume that the essential range of
$\Noise$ is dense in $S$. Then $C_\Noise$ is one-to-one, self-adjoint and {\ntext in the trace class}, and we may define the unique positive and self-adjoint
power $C_\Noise^{t}$ for any $t\in \real .$ For $t\geq 0$ we denote $S^t:=C_\Noise^{t}S$; henceforth
 the spaces $S^{1/2}$ and $S^1$ in
(\ref{basicdiagram}) are well defined.

The space $Z=S^{1/2}=C_\Noise^{1/2}S$ is called the {\em Cameron-Martin space} of $\Noise$, and $\Noise$ is called the white noise process
in $Z$. We remark that the realizations of $\Noise$ belong to $Z$ only with probability zero. Note that $C_\Noise^{t}:S\to S^t$ is an
isomorphism when the norm of the dense subspace $S^t\subset S$ is defined as $\| C_\Noise^{t}u\|_{S^t}:=\| u\|_S.$ The domain of definition
of $C_\Noise^{-t}$ is $S^t,$ and the map $C_\Noise^{-t}:S^t\to S$ is an isomorphism. Then
\begin{eqnarray*}
  \bra z_1,z_2\cet_Z=\bra C^{-1/2}_\Noise z_1,C^{-1/2}_\Noise z_2\cet_S.
\end{eqnarray*}
In the Gaussian example in Section \ref{sec:Gaussianexample} we have $S^t=H^{-2+4t}(\torus^2)$ so that $S=S^0=H^{-2}(\torus^2)$ and
$Z=L^2(\torus^2)$ and we may choose $Y=H^{-1}(\torus^2)$. 
Considering $S^{t+\frac 12}\subset Z\subset S^{-t+\frac 12}$ as a Gel'fand triple,
where  $S^{t+\frac 12}$ and $S^{-t+\frac 12}$ are considered as dual spaces
with the pairing $\bra w,z\cet_{S^{t+1/2}\times S^{-t+1/2}}=
\bra C_\Noise^{-t}w,C_\Noise^{t}z\cet_Z$,
it then holds that
\ba
\expec (\bra \Noise,z_1\cet_{S^0\times S^1} \bra \Noise,z_2\cet_{S^0\times S^1})= \bra
z_1,z_2\cet_Z\quad \hbox{for }z_1,z_2\in S^1.
\ea

Finally, we assume that $A:Y\to S^1$ is a bounded linear operator. The definition of {\ntext the} objects in diagram (\ref{basicdiagram}) is now
complete, and we turn to discussing reconstructors.

We analyze finite- and infinite-dimensional Bayesian inversion simultaneously, so let us introduce a rigorous setting for discrete
approximations to the random variable $U$ above.
We say that $Y$-valued random variables $U_n$ tend {\it weakly in distribution (w.i.d.)}
to $U$ if
\ba
\lim_{n\to\infty} \langle U_n ,y'\rangle =\langle U ,y'\rangle\quad \hbox{ in distribution}
\ea
for every $y'\in Y'$. {\nnntext Note that here ``weakly'' refers to the weak topology used in
space $Y$, and ``distribution'' to the convergence of scalar-valued  random variables.}

\begin{definition}\label{def:linear_discr}{\bf (Linear discretization of random functions)} Let $Y$ be a separable Banach space and $Y_n\subset Y$ be
finite-dimensional subspaces. {\mntext The spaces $Y_n$ need not be
nested.}
Let $U_n=U_n(\omega)$ be $Y_n$-valued random variables with $n=1,2,3,\dots$. Assume that
\begin{enumerate}
\item There is a $Y$-valued random variable $U=U(\omega)$ such that
\begin{eqnarray*}
  \lim_{n\to \infty} U_{n} = U \qquad \mbox{ w.i.d.}
\end{eqnarray*}

\item There are bounded linear {\nnntext operators}  $T_{n}\in L(Y)$ such that
$$
  U_{n}= T_{n} U.
$$
\end{enumerate}
Then we say that the $U_{n}$, $n\geq 1$, are {\em proper linear discretizations} of $U$ {\mntext in $Y$}.
\end{definition}

{\mntext Examples of random variables that form or do not form
proper linear discretizations are given in Appendix
\ref{sec:discretizationexamples}.}

{\nnntext In following, we mainly consider cases where $T_n$ are projection operators.}

Recall the following definition from \cite[II.7]{shir}.
\begin{definition}
Let $U\in L^1(\Omega,\Sigma;Y)$ be $Y$-valued random variable and $\Sigma_0$ a sub $\sigma$-algebra of $\Sigma$. Then the
conditional expectation $\expec(U|\Sigma_0)$ of $U$ exists with respect to $\Sigma_0$. That is,
 $\expec(U|\Sigma_0)\in L^1(\Omega,\Sigma_0;Y)$ and it satisfies
\begin{eqnarray*}
  \int_D \expec(U|\Sigma_0)(\omega)\,\prob(d\omega)=
  \int_D U(\omega)\,\prob(d\omega)\quad\hbox{for all $D\in \Sigma_0$} .
\end{eqnarray*}
\end{definition}
Note that all vector-valued integrals in this work are standard Bochner integrals. We refer the reader to
 \cite{du} for definition and basic facts on Bochner integral and vector-valued conditional expectations. The operator $P:U\mapsto \expec(U|\Sigma_0)$
is a projection $P:L^1(\Omega,\Sigma;Y)\to L^1(\Omega,\Sigma_0;Y)$, {where $L^1(\Omega,\Sigma_0;Y)$ denotes the space of
measurable functions $V:(\Omega,\Sigma_0)\to (Y,\tau^n)$ which are Bochner integrable.}

Now we are ready to give the definition of a (non-unique) reconstructor.

\begin{definition}\label{def: conditional mean function}
Denote by $\M\subset \Sigma$ be the $\sigma$-algebra generated by the random variable $M(\omega)$. We say that any deterministic function
$$
  \mathcal{R}_M(U|\,\,\cdotp\,):S\to Y,\qquad \qquad m\mapsto \mathcal{R}_M(U|\, m),
$$
is a {\em reconstructor} of $U$ with measurement $M$ if
$$
  \mathcal{R}_M(U|M(\omega)) = \expec(U|\M)(\omega)\quad \hbox{almost surely.}
$$
If $\widetilde Y$ is a separable Banach space and $g:Y\to \widetilde Y$ is a measurable function, we define
$\mathcal{R}_M(g(U)|\,\cdotp\,):S\to \widetilde Y$ to be any deterministic function satisfying
$$
   \mathcal{R}_M(g(U)|M(\omega))=\expec(g(U)|\M)(\omega)\quad \hbox{almost surely.}
$$
\end{definition}

{\ntext Note that reconstructor is a deterministic function. Also, if
a realization of the measurement in the computational model, $M_{kn}(\omega_0)$, is substituted
in the reconstructor ${\mathcal R}_M(U|\,\,\cdotp\,)$, then the obtained result
$\mathcal{R}_M(U|\, M_{kn}(\omega_0))$ coincides with conditional expectation. Thus,
reconstructor, {\nnntext considered as a deterministic function, is just the}
functional representation of conditional expectation  (see e.g. \cite{shir}, section II.7, formula (43)).
However, we assume that we are not given a realization of the measurement
in computational model, $M_{kn}(\omega_0)$, as data, but instead $M_k(\omega_0)=P_k M(\omega_0)$,
a realization of the measurement in the practical measurement model (\ref{eq: measmodel_k k}). An essential feature of the reconstructor is that
it is defined for all elements in $S$. Thus, even though
the reconstructor is related to the computational model (\ref{eq: measmodel_eka_kn}), it is possible to
substitute into the reconstructor
the realization of the practical measurement model (\ref{eq: measmodel_k k}).}
{\nnntext
This is the  reason why the reconstructor, which has the same functional representation as conditional expectation,
is defined as a new concept.}

One may generalize the well-known scalar-valued result on the existence of reconstructor, see
\cite[II.3 Theorem 3 and II.7.5]{shir}, to Bochner-valued conditional expectations. However, we will not
need this since we next
 establish a specific formula for a reconstructor in our situation. {\ntext Note that the following result
is close to the usual functional representation of the conditional expectation, c.f.\ \cite{shir}.
As we need to define
reconstructors as a deterministic function of the space $S$ and also introduce
notations for later use,  we present the proof of the result for completeness. }

\begin{theorem} \label{g-theorem}
Denote by $\mu:\B_S\to [0,1]$ the distribution of $\Noise$ in $S$, and set $\mu_a(E)=\mu(E- a)$ for $a\in S$, $E\in \B_S$. Let $g:Y\to Y$
be a measurable function satisfying $\expec \|g(U)\|_Y <\infty.$
{\mntext Set for $m \in S$
\begin{eqnarray}\label{RN-formulae}
  H_{{U,M}}^g(m) &=&
  \int_{Y} g(u) \exp(-\frac 12 \| Au\|_Z^2+\bra C^{-1}_\Noise Au,m\cet_{S})\, d\lambda (u) ,\\ \nonumber
  H_{{U,M}}^1(m) &=&
  \int_{Y} \exp(-\frac 12 \| Au\|_Z^2+\bra C^{-1}_\Noise Au,m\cet_{S})\, d\lambda (u),
\end{eqnarray}
where $\lambda$ stands for the distribution of $U$ in $Y$.}
Then the function
\begin{eqnarray}\label{Matti eq: B}
  \mathcal{R}_M(g(U)| m)=\frac {H_{{U,M}}^g(m)} {H_{{U,M}}^1(m)},\quad m\in S,
\end{eqnarray}
is well-defined and satisfies $\mathcal{R}_M(g(U)|M(\omega))=\expec(g(U)|\M)(\omega)$ almost surely,
that is, formula (\ref{Matti eq: B}) defines
a reconstructor.
\end{theorem}

\begin{proof} 
Using the equality $M=AU+\Noise$, and Fubini theorem, 
 we have for any measurable $D=E\times F\subset S\times
Y$
\begin{eqnarray}\label{Z}
  \prob(\{(M,U)\in D\})
  &=&\int_{\Omega_1}\left[ \int_S  \chi_{{}_E}(\noise+AU(\omega_1)) \, \chi_{{}_F}(U(\omega_1))\,
  d\mu(\noise)\right]d\prob_1(\omega_1)\nonumber\\
  &=&\int_{\Omega_1}\left[ \int_S  \chi_{{}_E}(m') \, \chi_{{}_F}(U(\omega_1))\, d\mu_{AU(\omega_1)}(m')\right]d\prob_1(\omega_1).\nonumber
\end{eqnarray}
{\ntext Above, $m'$ is an integration variable running over the space $S$ where $M$ is taking values}.
 Thus  we have for any integrable function $f:S\times
Y\to \C$
\begin{eqnarray}\label{Matti eq: A}
\expec(f(M,U))=\int_{\Omega_1}\left[ \int_S f(m',U(\omega_1)) d\mu_{AU(\omega_1)}(m')\right]d\prob_1(\omega_1).
\end{eqnarray}
One checks that the same holds for any Bochner integrable function $f:S\times Y\to Y$ by simply using the fact that such an $f$ is an almost sure  limit of simple functions $f_k$ that satisfy the pointwise inequality
 $\| f_k\|_Y\leq \| f \|_Y$.

Since $Z$ is the Cameron-Martin space of $\Noise$, we have for any $a\in S^1$ by \cite[Cor. 2.4.3]{B}
the  Radon-Nikodym\footnote {To
motivate formula (\ref{translation a}), we note that if  $Z$ would be finite-dimensional, we could write
\begin{eqnarray*}
 \frac {d \mu_a} {d \mu}(m)=
\exp(-\frac 12 \|m-a\|_Z^2 + \frac 12 \|m\|_Z^2) = \exp(-\frac 12 \|a\|_Z^2+\bra a,m\cet_{Z}).
\end{eqnarray*}
The formula (\ref{translation a}) is a generalization of this
 for the infinite-dimensional case.} derivative
\begin{eqnarray}\label{translation a}
 \frac {d \mu_a} {d \mu}(m)&=&
\exp(-\frac 12 \|a\|_Z^2+\bra C^{-1}_\Noise a,m\cet_{S}).
\end{eqnarray}
The latter formula has the advantage of being well-defined for every $m\in S.$ In particular, we have
\begin{eqnarray}\label{translation}
 \frac {d \mu_{AU(\omega_1)}} {d \mu}(m)=
\exp\Big(-\frac 12 \|AU(\omega_1)\|_Z^2+\bra C^{-1}_\Noise AU(\omega_1), m\cet_{S}\Big).
\end{eqnarray}

{\mntext Using formula (\ref{translation})
we see that formula (\ref{RN-formulae})
can be written as
\begin{eqnarray*}
  H_{{U,M}}^g(m) &=&  \int_{\Omega_1} g(U(\omega_1)) \frac {d \mu_{AU(\omega_1)}} {d \mu}(m)
  \,d\prob_1(\omega_1),\\
  H_{{U,M}}^1(m) &=&  \int_{\Omega_1} \frac {d \mu_{AU(\omega_1)}} {d \mu}(m)
  \,d\prob_1(\omega_1).
\end{eqnarray*}
}

By Fubini theorem we may continue from \refeq{Z} to obtain
\begin{eqnarray*}
  \expec(f(M,U))=\int_S\left[ \int_{\Omega_1} f(m',U(\omega_1)) \frac {d \mu_{AU(\omega_1)}} {d \mu}(m') \,d\prob_1(\omega_1) \right]  d\mu(m').
\end{eqnarray*}
Especially for $E\subset \B_S$ it holds that
\begin{eqnarray*}
  \expec(\chi_E(M)g(U))&=&\int_E\left[ \int_{\Omega_1} g(U(\omega_1)) \frac {d \mu_{AU(\omega_1)}} {d \mu}(m') \,d\prob_1(\omega_1) \right]
  d\mu(m').
\end{eqnarray*}
Let $\nu:\B_S\to [0,1]$ be the measure $\nu(E)=\prob(M^{-1}(E))$, that is, $\nu$ is the distribution of $M$. Now
\begin{eqnarray*}
  \int_E d\nu(m')=\expec(\chi_E(M))&=&\int_E\left[ \int_{\Omega_1} \frac {d \mu_{AU(\omega_1)}} {d \mu}(m') \,d\prob_1(\omega_1) \right]
  d\mu(m')
\end{eqnarray*}
and thus $\nu <\!< \mu$ and $\frac {d\nu}{d\mu}(m)=H_{{U,M}}^1(m)$ for almost every $m$ with respect to $\mu$, where
\begin{eqnarray*}
  H_{{U,M}}^1(m)=\int_{\Omega_1} \frac {d \mu_{AU(\omega_1)}} {d \mu}(m) \,d\prob_1(\omega_1).
\end{eqnarray*}
Observe that $H^1_{U,M}(m)>0$ for all $m\in S$ and thus $\nu(E)=0$ if and only if $\mu(E)=0$. Hence also $\mu \ll\nu$ and
\begin{eqnarray*}
  \frac {d\mu}{d\nu}(m)=\Big(\frac {d\nu}{d\mu}(m)\Big)^{-1}
\end{eqnarray*}
is well defined. Now  by (\ref{Matti eq: A})
\begin{eqnarray*}
  \int_E \int_{\Omega_1} \|g(U(\omega_1))\|_Y \frac {d \mu_{AU(\omega_1)}} {d \mu}(m') (\frac {d\nu}{d\mu}(m'))^{-1}\,d\prob_1(\omega_1)
  d\nu(m')
  \leq\expec(\|g(U)\|_Y)<\infty.
\end{eqnarray*}
Hence by Fubini theorem for $E\in \B_S$
\begin{eqnarray*}
  & &\expec(\chi_E(M)g(U))\\
  &=&\int_E\left[ \int_{\Omega_1} g(U(\omega_1)) \frac {d \mu_{AU(\omega_1)}} {d \mu}(m') \,d\prob_1(\omega_1) \right](\frac
  {d\nu}{d\mu}(m'))^{-1}  d\nu(m')
  \\
  &=&\int_\Omega \chi_E(M(\omega))\left[ \int_{\Omega_1} g(U(\omega_1)) \frac {d \mu_{AU(\omega_1)}} {d \mu}(M(\omega))
  \,d\prob_1(\omega_1)
  \right](\frac {d\nu}{d\mu}(M(\omega)))^{-1}  d\prob(\omega).
\end{eqnarray*}
Thus we have the almost sure equality
\begin{eqnarray*}
  \expec(g(U)|\M)(\omega)= \left[ \int_{\Omega_1} g(U(\omega_1)) \frac {d \mu_{AU(\omega_1)}} {d \mu}(M(\omega))
  \,d\prob_1(\omega_1)
  \right](\frac {d\nu}{d\mu}(M(\omega)))^{-1}.
\end{eqnarray*}
This verifies the formula for  $\mathcal{R}_M(g(U)|m)$ given in the assertion.
 \hfill$\Box$
\medskip
\end{proof}

For convenience, let us look at Theorem \ref{g-theorem} in the finite-dimensional case. Then $Y=\R^n$ and $Z=S=\R^m$, $\Noise$ is the
Gaussian white noise with the identity covariance matrix, and $U$ and $M$ have smooth everywhere positive probability density functions
$\pi_U(u)$ and $\pi_M(m)$, correspondingly. It follows that
\begin{eqnarray}\label{fin dim formula}
  \mathcal{R}_M(U| m)=(2\pi)^{-n/2}\int_{\R^n}
  u \exp(-\frac 12 \|Au-m\|_Z^2)\, \frac {\pi_U(u)}{\pi_M(m)}du
\end{eqnarray}
satisfies the assumptions of Definition \ref{def: conditional mean function}. Compare formulas (\ref{Matti eq: B}) and (\ref{fin dim
formula}) and see Appendix \ref{sec:recDomain} for further discussion. Note that (\ref{fin dim formula}) is widely used in practical Bayesian inversion \cite{tarantolabook,KS,EkiDanielaBook}.

Stability of the reconstructor with respect to data $m$ is important from the point of view of practical inversion. The following
theorem yields a non-quantitative convergence result for reconstructors in general. We provide sharper results for the special case of
Besov priors later in sections \ref{sec:technical} and \ref{sec:BesovConv}.

\begin{theorem} \label{convergence -theorem}
Assume that the exponential moments of $U$ satisfy
\begin{eqnarray}\label{abc}
  \expec(\exp(\lambda\| U\|_Y)) <\infty \quad\mbox{for all} \quad \lambda >0.
\end{eqnarray}
Take $g:Y\to \R$ to be such a continuous function that
\begin{eqnarray}\label{aaa}
|g(u)| \leq a\exp (a\|u\|_Y)\quad \mbox{for}\;\; u\in Y
\end{eqnarray}
with some constant $a$. Assume that $T_n:Y\to Y$, $n>0$, and $P_k:S^1\to S^1$, $k>0$, are linear projections satisfying
\begin{eqnarray}\label{eq:T limit}
  & &\lim_{n\to \infty}
  \|T_n y-y\|_{Y}=0\quad\hbox{ for all }y\in Y,\\
   \label{eq:P limit}
  & &\lim_{k\to \infty}
  \|P_k z-z\|_{S^1}=0\quad\hbox{ for all }z\in \hbox{Ran}(A),\\
\label{eq:C_0 bounds}
& &\|T_n\|_{L(Y)}\leq C_0,\quad \|P_k\|_{L(S^1)}\leq C_0\quad \hbox{for all $n,k$}\end{eqnarray}
with some $C_0>0$.
Finally, let $m_k,m\in S$ satisfy
\begin{eqnarray}\label{BBB formula}
  \lim_{k\to \infty} m_k=m\quad \hbox{in }S.
\end{eqnarray}
Then we have the convergence
\begin{eqnarray*}
  \lim_{n,k\to \infty}\mathcal{R}_{\Meas_{kn}}(g(U_n)| m_k ) =\mathcal{R}_M(g(U)|m),
\end{eqnarray*}
where the reconstructors are defined using formula (\ref {Matti eq: B}) for models (\ref{eq: measmodel_kn}) and (\ref{basicequation A1}),
respectively. Moreover,  the limits
$$
  \lim_{n\to \infty}\mathcal{R}_{\Meas_{kn}}(g(U_n)| m_k ) \quad\mbox{and}
  \quad   \lim_{k\to \infty}\mathcal{R}_{\Meas_{kn}}(g(U_n)| m_k )
$$
exist for a fixed value of $k$ (resp. $n$).
\end{theorem}

\begin{proof}
We have
\begin{eqnarray*}
  H_{U_n,\Meas_{kn}}^1(m_k)=\int_{\Omega_1} \frac {d \mu_{ P_kAT_n U(\omega_1)}} {d \mu}(m_k) \,d\prob_1(\omega_1)
\end{eqnarray*}
and  since \refeq{abc} and \refeq{aaa}  imply $\expec |g(T_nU)| <\infty$ we may also write
\begin{eqnarray*}
  H_{U_n,\Meas_{kn}}^g(m_k)=\int_{\Omega_1}g(T_nU(\omega_1)) \frac {d \mu_{ P_k A T_n U(\omega_1)}} {d \mu}(m_k) \,d\prob_1(\omega_1).
\end{eqnarray*}
Above
\begin{eqnarray}\label{Y}
  \frac {d \mu_{P_kAT_n U(\omega_1)}} {d \mu}(m_k)&=& \exp(-\frac 12 \|P_kAT_n U(\omega_1)\|_Z^2+\bra C^{-1}_\noise P_kAT_n U(\omega_1), m_k\cet_{S})\nonumber\\
&\leq & \exp (c\|U(\omega_1)\|_Y)\nonumber \quad \mbox{for all}\;\; k,n,
\end{eqnarray}
by our assumptions. The claims follows now from the Lebesgue dominated convergence theorem by applying the majorant $a\exp
((c+aC_0)\|U(\omega_1)\|_Y).$
 \hfill$\Box$
\medskip
\end{proof}

Now the general theory of discretization-invariant Bayesian inversion is in place for the case of measurement models $M = AU + \Noise$ and $\Meas_{kn} = P_kAT_nU + \Noise$ concerning infinite-dimensional noise $\Noise$. Using the continuum noise $\Noise$ is convenient above because we can work in the same function space regardless of $k$.

{\ntext Assume given data $M_k(\omega_0)$ corresponding to the practical measurement model (\ref{eq: measmodel_k k}) and consider the computational solution of the inversion problem using the computational model (\ref{eq: measmodel_eka_kn}),} where random error is finite-dimensional white noise $\Noise_k=P_k\Noise$.  It remains to discuss the implications of our general theory for these practical models. To do this, assume that $P_k$ are projections $P_k:S\to S$ having the following
properties:
\beq
  &&\label{Pk 1}\hbox{Ran}\,(P_k)\hbox{ is a finite-dimensional subset of }S^1,\\
  &&\label{Pk 2}\bra P_k \phi,\psi\cet_Z=\bra  \phi,P_k\psi\cet_Z\quad\hbox{for }\phi,\psi\in Z.
\eeq
First we show that reconstructors corresponding to the measurement models $\Theta_{kn} = P_k AU_n+\Noise$ and $ M_{kn} = P_kAU_n+P_k \Noise$ actually coincide.
\begin{lemma}\label{reconstructors coincide}
Assume $P_k:S\to S$ is a projection satisfying (\ref{Pk 1}) and (\ref{Pk 2}). Then the reconstructors  defined in Theorem \ref{g-theorem}
for the  measurement models $\Theta_{kn} = P_k AU_n+\Noise$ and $ M_{kn} = P_kAU_n+P_k \Noise$ satisfy
\ba
  \mathcal{R}_{\Theta_{kn}}(g(U_n)|m_k)=
  \mathcal{R}_{M_{kn}}(g(U_n)|m_k)
\ea
for $m_k\in \hbox{Ran}\,(P_k)$.
\end{lemma}

\begin{proof}
Consider first $\Noise_k=P_k\Noise$ as a Gaussian random variable
 taking values in the space
$\hbox{Ran}\,(P_k)$ that has the inner product inherited from $S$.
The random variable $\Noise_k$ has zero expectation.

Consider $S$ and $S^1$ as dual Hilbert spaces. The corresponding pairing is
\ba
\bra  \eta,\phi\cet_{S\times S^1}
=\bra  \eta,C^{-1}_\Noise\phi\cet_S,\quad \eta\in S,\ \phi\in S^1,
\ea
which is an extension of the pairing
\ba
\bra \eta,\phi\cet_{S\times S^1}=\bra \eta,\phi\cet_Z
\ea
defined for $\eta\in Z\subset S$ and $\phi\in S^1$. Moreover,
for such $\eta$ and $\phi$
\ba
\bra P_k\eta,\phi\cet_{S\times S^1}=\bra P_k\eta,\phi\cet_Z
=\bra \eta,P_k\phi\cet_Z=\bra \eta, P_k\phi\cet_{S\times S^1}.
\ea
 As the finite-dimensional projection $P_k:S\to S$
is bounded, the density of $Z$ in $S$ implies that
$\bra P_k\eta,\phi\cet_{S\times S^1}=\bra \eta,P_k\phi\cet_{S\times S^1}$
for all $\eta\in S$ and $\phi\in S^1$.
Using this, we see for $\phi,\psi\in \hbox{Ran}\,(P_k)$
\ba
\expec (\bra \Noise_k,\phi\cet_Z\bra \Noise_k,\psi\cet_Z)&=&
\expec (\bra P_k\Noise,\phi\cet_{S\times S^1}\\
&=&
\bra P_k\Noise,\psi\cet_{S\times S^1})\\
&=&
\expec (\bra \Noise,P_k\phi\cet_{S\times S^1}\bra\Noise,
P_k\psi\cet_{S\times S^1})\\
&=&
\expec (\bra \Noise,C_\Noise^{-1}P_k\phi\cet_S\bra\Noise,
C_\Noise^{-1}P_k\psi\cet_S)\\
&=&\bra C_\Noise^{-1} \phi,
C_\Noise C_\Noise^{-1}\psi\cet_S\\
&=&\bra \phi,\psi\cet_Z.
\ea
This implies that the covariance operator of $\Noise_k$,
considered now as a Gaussian random variable taking values in $\hbox{Ran}\,(P_k)$
endowed with the inner product inherited from $Z$, is
the identity operator.
Using this we see that the reconstructor  defined in Theorem \ref{g-theorem}
for the measurement $M_{kn}$ has the form
\ba
  \mathcal{R}_{M_{kn}}(g(U_n)|m_k)
&=&
  \frac{\expec(g(U) \exp(-\frac 12 \| AU_k-m_k\|_Z^2))}
  {\expec( \exp(-\frac 12 \| AU_k-m_k\|_Z^2))}
\\
&=&
  \frac{\expec(g(U) \exp(-\frac 12 \| AU_k\|_Z^2+\bra A U_k,m_k\cet_{Z}))}
  {\expec( \exp(-\frac 12 \| AU_k\|_Z^2+\bra A U_k,m_k\cet_{Z}))}
\ea
for $m_k\in \hbox{Ran}\,(P_k).$
The reconstructor $\mathcal{R}_{\Meas_{kn}}(g(U_n)|m_k)$
defined in Theorem \ref{g-theorem} for the measurement $\Meas_{kn}$
has the same form, and the assertion follows.
 \hfill$\Box$
\medskip
\end{proof}

Finally, we prove the convergence of reconstructors for models with discrete noise.

\begin{theorem}\label{Corollary-teoreema}
Assume that in addition to conditions (\ref{abc})-(\ref{eq:C_0 bounds}),
the projections $P_k:S\to S$ satisfy (\ref{Pk 1}), (\ref{Pk 2}),
together with
\beq
 \label{eq:P limit 2}
 \lim_{k\to \infty}
  \|P_k z-z\|_{S}=0\quad\hbox{ for all }z\in S.
\eeq
Let $u=U(\omega_0)$,  $\e=\Noise(\omega_0)$, $\omega_0\in \Omega$ be realizations of the random variables $U$ and $\Noise$, and let
\ba
 m = A u + \e,\qquad
 m_k = A_k u+P_k\e,
\ea
be the realizations of the measurements (\ref{basicequation A1}) and (\ref{eq: measmodel_k k}), respectively. Then
the reconstructors  defined in Theorem \ref{g-theorem}
for the measurement models $M_{kn} = P_kAT_nU + P_k\Noise$ and $M = AU + \Noise$
 satisfy
\ba
 & & \lim_{n,k\to \infty}\mathcal{R}_{M_{kn}}(g(U_n)| m_k ) =\mathcal{R}_M(g(U)|m).
\ea
Moreover,  the limits
$$
  \lim_{n\to \infty}\mathcal{R}_{M_{kn}}(g(U_n)| m_k ) \quad\mbox{and}
  \quad   \lim_{k\to \infty}\mathcal{R}_{M_{kn}}(g(U_n)| m_k )
$$
exist for a fixed value of $k$ (resp. $n$).
\end{theorem}
\begin{proof}
By (\ref{eq:P limit 2}), $\lim_{k\to \infty} m_k=m$ in $S$,
and hence the assertion follows by Theorem \ref{convergence -theorem}
and Lemma \ref{reconstructors coincide}.
 \hfill$\Box$
\end{proof}

\medskip
Theorem \ref{Corollary-teoreema} concerns the convergence of practical inversion methods: $m_k$ is data provided by an actual measurement device, and the computational model $M_{kn} = P_kAT_nU + \Noise_k$ allows computer implementation. For instance, most Markov chain Monte Carlo inversion algorithms are programmed to evaluate $\mathcal{R}_{M_{kn}}(U_n| m_k )$.
\medskip

Let $E\subset Y$ be a Borel set
and $\chi_{{}_E}$ be the indicator function of $E$. Using 
reconstructors, we define 
\ba
& &{\mathcal P}(E|m)=\mathcal{R}_{M}(\chi_{{}_E}(U)| m),\\
& &{\mathcal P}_{kn}(E|m_k)=\mathcal{R}_{M_{kn}}(\chi_{{}_E}(U_n)| m_k).
\ea
For a given $m\in S$,
the map $E\mapsto  {\mathcal P}(E|m)$ is a probability measure on $Y$
by equation (\ref{Matti eq: B}). Next, let $E$ be a fixed Borel set of $Y$.
If we substitute $M(\omega)$ in the function
$m\mapsto {\mathcal P}(E|m)$, we obtain by Definition \ref{def: conditional mean function}
\ba
{\mathcal P}(E|M(\omega))=\expec(\chi_{{}_E}(U)|\M)(\omega)
=\prob(\{U\in E\}|\M)(\omega)\quad\hbox{almost surely},
\ea
where $\prob(\{U\in E\}|\M)(\omega)$ is 
the conditional probability for the event 
$\{U\in E\}$ with respect to the $\sigma$-algebra $\M$
(see \cite[Sec.\ 6]{kall}).
Thus, roughly speaking, $m\mapsto {\mathcal P}(E|m)$ can be considered as 
the posterior probability for the event $U\in E$ when
the measurement $M$ gets value $m$, and
${\mathcal P}_{kn}(E|m_k)$ the posterior probability for the event $U_n\in E$ when
the measurement $M_{kn}$ gets value $m_k$.

Recall that  measures $\mu_j$ on $Y$ convergence weakly
to measure $\mu$ as $j\to \infty$
if $\lim_{j\to \infty}\int_Y g\,d\mu_j=\int_Y g\,d\mu$ for all 
bounded continuous functions $g:Y\to \R$. Thus 
Theorem \ref{Corollary-teoreema} yields the  following corollary.

\begin{corollary}\label{co:10}  Let the assumptions
of Theorem \ref{Corollary-teoreema} hold. Then the Borel measures
$E\mapsto {\mathcal P}_{kn}(E|m_k)$ converge
weakly to the measure $E\mapsto {\mathcal P}(E|m)$
as $n,k\to \infty$.
\end{corollary}
\medskip

\section{Random variables in Besov spaces}\label{sec:randomBesov}

We wish to use priors of the form
$$
  \prior(u)\underset{\mbox{\tiny \em formally}}= c\exp(-\| u\|_{B_{pp}^s(\torus^d)}^p)
$$
for any integrability parameter $1\leq p<\infty$ and some chosen smoothness $s\in\R$. {\mntext Note that we consider now functions on a general $d$-dimensional torus $\T^d$.} 

Following Appendix \ref{app:wavelet} we use a  compactly supported wavelet $\tilde{\psi}(x)$, $x\in \R^d$
and scaling function $\tilde{\phi}(x)$ suitable for
multi-resolution analysis of smoothness $C^r$ in $L^2(\R^d)$ with large enough $r$. Then we can expand functions as
$$
  f(x) = \sum_{\ell=1}^\infty c_\ell \psi_\ell(x),\quad x\in \T^d
$$
and $f\in B^{s}_{pp}(\T^d)$ if the norm \cite{meyer}
\begin{equation}\label{normi1}
  \|\sum_{\ell=1}^\infty c_\ell \psi_\ell(x)\|_{B^{s}_{pp}(\T^d)} :=\left(\sum_{\ell=1}^\infty \ell^{(ps/d+p/2-1)}|c_\ell|^p \right)^{1/p}
\end{equation}
is finite. For the exact bookkeeping of scales and locations of $\psi_{\ell}$ as function of $\ell$ we refer to Appendix \ref{app:wavelet}.

\begin{definition}\label{besov-Gaussian}
Let $1\leq p<\infty$ and $s\in\R$. Let $(X_\ell)_{\ell=1}^\infty$ be  independent identically distributed
{\ntext real-valued random variables with probability density function}
\begin{eqnarray}\label{X-distri}
  \pi_X(x)=c_p\exp(-|x|^p),\quad \hbox{\ntext with $c_p=\left(\int_\R\exp(-|x|^p)\,dx \right)^{-1}$}.
\end{eqnarray}
Let $U$ be the random function
\begin{eqnarray*}
  U(x)=\sum_{\ell= 1}^\infty \ell^{-(s/d+1/2-1/p)}X_{\ell}\psi_{\ell}(x),\quad x\in \torus^d .
\end{eqnarray*}
Then we say that $U$ is distributed according to {\ntext a} $B_{pp}^s$ prior.
\end{definition}

Next we show that random variable $U$ in Definition \ref{besov-Gaussian} is a well defined object.

\begin{lemma}\label{le1} Let $U$ be as in Definition \ref{besov-Gaussian} and take $t\in\R$. The following three conditions are equivalent:

\smallskip

\noindent {\rm (i)} \qquad $\| U\|_{B_{pp}^{t }} <\infty$ \quad \mbox{almost surely}.

\smallskip

\noindent {\rm (ii)} \qquad $\ds \expec \exp \big(\frac 12  \| U\|_{B_{pp}^{t}}^p\big)<\infty .$

\smallskip

\noindent {\rm (iii)} \qquad $t <s-\frac dp$.

\end{lemma}
\begin{proof} Denote by $(X_\ell)_{\ell=1}^\infty$ a sequence of {\ntext independent random
variables} with density \refeq{X-distri}. Assume (iii). In order to deduce (ii) we need to show the finiteness of the quantity
\begin{eqnarray}
  \expec\exp\Bigl( {1\over 2}\|\sum_{\ell =1}^\infty \ell^{-s/d-(1/2-1/p)}X_\ell\psi_\ell\|_{B_{pp}^{t }}^p\Bigr)
  &\simeq& \nonumber
  \expec\exp\Bigl(\sum_{\ell =1}^\infty \frac 12 \ell^{(t -s)p/d}|X_\ell |^{p}\Bigr)\\
  &=& \nonumber
  \prod_{\ell =1}^\infty\expec \exp (\frac 12 \ell^{-(s-t )p/d}|X_\ell |^{p}) \\
  &=& \label{eq106}
  \prod_{\ell =1}^\infty(1-\ell^{-(s-t )p/d}/2)^{-1/p}.
\end{eqnarray}
Above we used independence and  the observation that for $k\in (0,1)$ one may compute $\expec \exp\big(k|X_\ell |^{p}\big)=(1-k)^{-1/p}$.
Clearly the product in (\ref{eq106}) converges if $t <s-d/p.$ The notation $a\simeq b$ stands for the existence of a positive constant $c<\infty$
such that $a/c\leq b\leq ca .$

Observe next that obviously (ii) implies (i). Finally, assume that (i) is true. Then
$$
\sum_{\ell =1}^\infty \ell^{(t -s)p/d}|X_\ell |^{p}<\infty
$$
almost surely. Since the the random variables $|X_\ell |^{p}$ are non-negative and identically distributed, an
easy application of {\ntext truncation} and \cite[Thm.\ 4.17]{kall}
 shows that almost sure finiteness of the sum  implies finiteness of the expectation. Hence (i)
implies (iii). \hfill$\Box$
\medskip
\end{proof}

Now we easily see that $B_{pp}^s(\torus^d)$ distributions generate discretization-invariant priors.

\begin{proposition} \label{discret inv}
Let $U$ be distributed according to $B_{pp}^s$ prior as in Definition \ref{besov-Gaussian}, and let $t<s - \frac dp$. Take $T_n:Y\to Y$
to be the sequence of linear operators defined by
$$
  T_n\left(\sum_{\ell=1}^\infty c_\ell \psi_\ell\right) =  \sum_{\ell=1}^n c_\ell \psi_\ell;
$$
then $\lim_{n\to \infty}T_n y=y$ for all $y\in Y:=B_{pp}^{t}(\T^d)$. Then $U_n=T_nU$, $n=1,2,\dots$ are proper linear discretizations of
$U$ in the sense of Definition {\rm \ref{def:linear_discr}}.
\end{proposition}

\begin{proof}
The random variables $U_n=T_nU$ converge almost surely in norm topology of $B_{pp}^{t}(\T^d)$. Since almost {\ntext sure} convergence implies convergence
weakly in distribution, the assertion follows. \hfill$\Box$
\medskip
\end{proof}
Note that if above $t>\frac dp$, then the continuous embedding $B_{pp}^{t}(\T^d)\to C(\T^d)$ implies that realizations of $U$ and $U_n$
are almost surely continuous.

\section{Quantitative estimates for reconstructors}\label{sec:technical}

This section studies the case that $U$ is distributed according to $B_{pp}^s$-prior.
We provide quantitative stability estimates for reconstructors. For $p>1$ qualitative results are
described by Theorem \ref{convergence -theorem}, thanks to Lemma \ref{le1}(ii). However, the case $p=1$ is more difficult  since condition \refeq{abc} fails.

We collect a set of assumptions together for later reference:
\medskip

\noindent{\bf Assumption A.}\quad{\sl Let $s\in\real$ and $1\leq p<\infty$ be arbitrary. Fix $t,\widetilde t,r\in\R $ such
that
$$
 t<\widetilde t <s-d/p\quad\mbox{ and}\quad
r >d/2.
$$
Let $A$ be a linear operator satisfying
\begin{eqnarray}\label{eq: basic assumption on A}
  A: B_{pp}^{t}(\torus^d)\to B_{11}^{r }(\torus^d). 
\end{eqnarray}
Assume that $g:B_{pp}^{t}\to B_{pp}^{t}$ is a map such that for some $q<\infty$ we have
\begin{eqnarray}\label{1}
&&\| g(u)\|_{B_{pp}^{t}}\leq c(1+ \| u\|_{B_{pp}^{t}})^q\quad\mbox{and}\nonumber\\
&&\| g(u_1)-g(u_2 )\|_{B_{pp}^{t}}\leq c\|u_1-u_2\|_{B_{pp}^{t}}(1+ \max(\| u_1\|_{B_{pp}^{t}},\| u_2 \|_{B_{pp}^{t}}))^q.\nonumber
\end{eqnarray}
 Finally, let $G: B^{\widetilde t}_{pp}\to B^{{t}}_{pp}$ be a bounded operator.}
\medskip

How does the diagram (\ref{basicdiagram}) look for the $B_{pp}^s$-prior? Take $t<s-\frac dp$ and set
$$
\begin{array}{cccccc}
&&&&&\\
 Y\!=\!B^{t}_{pp}(\T^d) &  \stackrel{A}{\longrightarrow}
& \ B^r_{11}(\T^d) & \ \subset\quad & B^{-r}_{\infty\infty}(\T^d).\\ \\
\rotatebox{90}{$\in$} &&&& \!\!\!\!\!\!\rotatebox{90}{$\in$}\\
U(\omega_1) &&&& \!\!\!\!\!\!\Noise (\omega_2)
\\ &&&&
\end{array}
$$
Observe that by Lemma \ref{le1} the random variable $U$ takes values 
in the Besov space $Y=B^{t}_{pp}(\T^d)$. Above $\Noise$ is standard 
Gaussian white noise: $\expec \Noise=0$ and 
\ba
\expec (\bra \Noise,\phi\cet\bra \Noise,\psi\cet)= \bra \phi,\psi\cet
\ea 
for all $\phi,\psi\in C^\infty(\T^d)$,
where $\bra\cdotp,\cdotp\cet$ is the distribution duality. 
To make our results more precise, we will consider the white noise
as a random variable taking values in a  
Besov space instead of Sobolev spaces, and 
 consider $\Noise$ here as a random variable taking values in the Besov space
$B^{-r}_{\infty\infty}(\T^d)$.

In assumption A we are particularly interested in $g$ being a characteristic function $g(u)=\chi_E(u)$ of some set $E$ (corresponding to confidence intervals), the identity map $g(u)=u$ (corresponding to the mean), and $g(u)=\|u\|_{B^{t}_{pp}}^q$.
We denote by $H^g(m , A, G)$ the quantity
$$
H^g( m , A, G):=\expec\left( g(GU)\exp \big(-{1\over 2}\| AU\|^2_{L^2}+\langle  AU,m\rangle_{} \big)\right).
$$
In the case $g\equiv 1$ the operator $G$ plays no role, and the corresponding real-valued quantity is denoted simply by $H^1(
m , A)$. In general, the above integral is defined as a $B_{pp}^t$-valued Bochner integral. The weak measurability of the
integrand is obvious, whence the strong measurability follows by the separability of the spaces involved. The following auxiliary result
will be used to verify integrability and to estimate the sensitivity of the quantity $H^g(m , A, G)$ with respect to changes in
the variables. Below, for $1\leq p< \infty$ we denote $p^\prime=p/(p-1)$.


\begin{proposition}\label{pr1}
Let $U$ be distributed according to $B_{pp}^s$ prior as in Definition \ref{besov-Gaussian}. Let $q\geq 0$ and
$m\in B^{-r}_{\infty\infty}$. Denote by $S_q(m , A)$ the random variable
$$
  S_q(m , A):=(1+\| U\|_{B^{\widetilde t }_{pp}})^{q}\exp \big(-{1\over 2}\| AU\|^2_{L^2}+\langle  AU,m\rangle_{} \big).
$$
Let $ w\geq r$ be an arbitrary index, which in case $p=1$ satisfies $ w > r$. 
Assume that $A:B^{t}_{pp}(\T^d)\to B^{w}_{11}(\T^d)$ is bounded.
Then there is a constant $c=c(p,r,t,w,d,q)$ such that
\begin{eqnarray}\label{eq100}\quad\quad
  \expec S_q(m , A)\leq c \exp \left( c( \|A\|_{B^{t}_{pp}\to B^{w}_{11}})^{2r/(w-r+2r/p^\prime)}(\| m
  \|_{B^{-r}_{\infty\infty}})^{{2w/(w-r+2r/p^\prime)}} \right)
\end{eqnarray}
with the understanding that in the case $p=1$ one sets $2r/p^\prime=0.$
\end{proposition}
\begin{proof}
In order to prove \refeq{eq100} we  first  observe that by Lemma \ref{le1}(ii) and {\ntext Cauchy-Schwarz} inequality it is enough to estimate the
expectation \beqla{iprime} I^\prime :=\expec\exp \big(-\| AU\|^2_{L^2}+2\langle  AU,m\rangle_{} \big). \eeq
 By Lemma \ref{le1} we have that $\expec\exp ({1\over 2}\| U\|^p_{B^t_{pp}})<\infty$.
The decomposition
$$
I^\prime =\expec\exp \big({1\over 2}\| U\|^p_{B^t_{pp}}+(-{1\over 2}\| U\|^p_{B^t_{pp}}-\| AU\|^2_{L^2}+2\langle  AU,m\rangle_{})
\big)
$$
shows that
\beq\label{EQ: A1}
I^\prime \leq c\,\exp\left({\sup_US(U)}\right),
\eeq
 where
$$
S(U):=-{1\over 2}\| U\|^p_{B^t_{pp}}-\| AU\|^2_{L^2}+2\langle  AU,m\rangle_{}.
$$
For $\ell\geq 0$ denote by $R_\ell$  the standard projection to the first $\ell$ coordinates in the wavelet basis, $\ell\geq 0$, with
$R_0=0$ and $Q_\ell=I-R_\ell.$ Fix an integer ${\ell_0}\geq 0$ and divide further above $m= R_{\ell_0}m
+Q_{\ell_0}m$.

We now assume that $p>1.$ By completing the square, and applying  Young's  inequality in the form $-x^p/2+2xy\leq c_py^{p^\prime}$ we
obtain
\begin{eqnarray}\label{0'}
S(U)&=& -{1\over 2}\| U\|^p_{B^t_{pp}} +2\langle AU,Q_{\ell_0}m
\rangle_{}-\| AU\|^2_{L^2}+2\langle  AU, R_{\ell_0} m
\rangle_{}\nonumber\\
&\leq& -{1\over 2}\| U\|^p_{B^t_{pp}} +2\| U\|_{B^t_{pp}} \| A^*Q_{\ell_0}m\|_{B^{-t}_{p^\prime p^\prime}}
-\| AU-R_{\ell_0}m\|^2_{L^2} + \| R_{\ell_0}m\|^2_{L^2} \nonumber\\
&\leq& c_p\| A^*Q_{\ell_0}m\|^{p^\prime}_{B^{-t}_{p^\prime p^\prime}}+\| R_{\ell_0}m\|^2_{L^2}.
\end{eqnarray}
If $m=0$ there is nothing to prove, so assume that $m\not=0$. Recall that $w\geq r$ and consider first the situation
where
\begin{eqnarray}\label{1'}
\|A^*\|^{p^\prime}_{B^{-w}_{\infty\infty}\to B^{-t}_{p^\prime p^\prime}}\geq \Big(\|m\|_{B^{-r}_{\infty\infty}}
\Big)^{2-p^\prime},
\end{eqnarray}
where $1/p+1/p^\prime=1$. In this case we choose
\ba
\ell_0&:=&\big[\left( \|A^*\|^{p^\prime}_{B^{-w}_{\infty\infty}\to B^{-t}_{p^\prime p^\prime}}
\|m\|^{p^\prime-2}_{B^{-r}_{\infty\infty}}\right)^{d/(2r+p^\prime(w-r))}\big]+1 \\
&\simeq& \left( \|A^*\|^{p^\prime}_{B^{-w}_{\infty\infty}\to
B^{-t}_{p^\prime p^\prime}}\|m\|^{p^\prime-2}_{B^{-r}_{\infty\infty}}\right)^{d/(2r+p^\prime(w-r))}.
\ea
We may estimate
\begin{eqnarray}\label{2'}
\| A^*Q_{\ell_0}m\|_{B^{-t}_{p^\prime p^\prime}}&\leq& \| A^*\|_{B^{-w}_{\infty\infty}\to B^{-t}_{p^\prime p^\prime}}
\|Q_{\ell_0}m\|_{B^{-w}_{\infty\infty}}\nonumber\\
&\leq& \ell_0^{(r-w)/d} \| A^*\|_{B^{-w}_{\infty\infty}\to B^{-t}_{p^\prime p^\prime}}\|m\|_{B^{-r}_{\infty\infty}}.
\end{eqnarray}
Moreover, one easily verifies that
\begin{eqnarray}\label{3'}
\|R_{\ell_0}m\|^2_{L^2}\leq (2\ell_0)^{2r/d}\|m\|^2_{B^{-r}_{\infty\infty}}.
\end{eqnarray}
By invoking the choice of $\ell_0$ and applying the auxiliary estimates \refeq{2'} and \refeq{3'} we compute from \refeq{0'} that
$$
S(U)\leq c( \|A\|_{B^{t}_{pp}\to B^{w}_{11}})^{2r/(w-r+2r/p^\prime)}(\|m\|_{B^{-r}_{\infty\infty}})^{
{2w/(w-r+2r/p^\prime)}}
$$
which finishes the proof in the situation \refeq{1'}.  Above, we also used the observation that $\| A^*\|_{B^{-w}_{\infty\infty}\to
B^{-t}_{p^\prime p^\prime}}=\|A\|_{B^{t}_{pp}\to B^{w}_{11}}$. On the other hand, if \refeq{1'} is not valid, that is,
\begin{eqnarray}\label{eq: opposite case}
\|A^*\|^{p^\prime}_{B^{-w}_{\infty\infty}\to B^{-t}_{p^\prime p^\prime}}< \Big(\|m\|_{B^{-r}_{\infty\infty}} \Big)^{2-p^\prime},
\end{eqnarray}
 we
apply the choice $\ell_0 =0$ in \refeq{0'}  to estimate
\begin{eqnarray}\label{4'}
S(U)&\leq& c_p \|A^*m\|^{p^\prime}_{B^{-t}_{p^\prime ,p^\prime}}\\
\nonumber
&\leq& c\| A^*\|^{p^\prime}_{ B^{-r}_{\infty\infty}\to
B^{-t}_{p^\prime p^\prime}} \|m\|^{p^\prime}_{B^{-r}_{\infty\infty}}\nonumber\\
&\leq& c(\|A\|_{B^{r}_{pp}\to B^{w}_{11}})^{2r/(w-r+2r/p^\prime)}(\|m\|_{B^{-r}_{\infty\infty}})^{{2w/(
w-r+2r/p^\prime)}}.\nonumber
\end{eqnarray}
In case $w>r$ the last  inequality above followed from inequality (\ref{eq: opposite case}).

Consider next the case $p=1.$  Assuming first that
\begin{eqnarray}\label{5'}
\| A^*\|_{ B^{-w}_{\infty\infty}\to B^{-t}_{\infty\infty}} \|m\|^{p^\prime}_{B^{-r}_{\infty\infty}}\geq 1/4
\end{eqnarray}
we utilize the fact that $w>r$ and choose
\ba
\ell_0&:=&  \big[(4\, \|m\|_{B^{-r}_{\infty\infty}}\|A^*\|_{B^{-w}_{\infty\infty}\to B^{-t}_{\infty\infty}})^{d/(w-r )}\big]+1\\ &\simeq&
(4\, \|m\|_{B^{-r}_{\infty\infty}}\|A^*\|_{B^{-w}_{\infty\infty}\to B^{-t}_{\infty\infty}})^{d/(w-r )}
\ea
and observe that then \refeq{2'} verifies that  $\| A^*Q_{\ell_0}m\|_{B^{-t}_{\infty\infty}}\leq 1/4,$ which in turn yields
\begin{eqnarray}\label{6'}
S(U)&=& -{1\over 2}\| U\|_{B^t_{11}} +2\langle AU,Q_{\ell_0}m\rangle_{}-\| AU\|^2_{L^2}+2\langle  AU,R_{\ell_0}m\rangle_{L^2}\nonumber\\
&\leq& -\| AU-R_{\ell_0}m\|^2_{L^2} + \| R_{\ell_0}m\|^2_{L^2} \\ &\leq&  \| R_{\ell_0}m\|^2_{L^2}
\\ & \leq& (2\ell_0)^{2r/d}\|m\|^2_{B^{-r}_{\infty\infty}}
\nonumber\\&\leq& c(4 \|A\|_{B^{t}_{11}\to B^{w}_{11}})^{2r/(w-r)}
(\|m\|_{B^{-r}_{\infty\infty}})^{{2w/(w-r)}},\nonumber
\end{eqnarray}
where we also applied the estimate \refeq{3'}. Finally, in case \refeq{5'} does not hold we may choose $\ell_0=0$ above and obtain that
$$
  S(U)\leq c( \|A\|_{B^{t}_{pp}\to B^{w}_{11}})^{2r/(w-r+2r/p^\prime)}(\| m
  \|_{B^{-r}_{\infty\infty}})^{{2w/( w-r+2r/p^\prime)}}
$$
for all $U$. Summarizing, we have shown that in all the above cases
\beq\label{EQ: A2}
\sup_U S(U)\leq
 c_1( \|A\|_{B^{t}_{pp}\to B^{w}_{11}})^{2r/(w-r+2r/p^\prime)}(\| m
  \|_{B^{-r}_{\infty\infty}})^{2w/(w-r+2r/p^\prime)},
\eeq
with some $c_1,c_2>0$, which finishes the proof of the Proposition.
 \hfill$\Box$
\medskip
\end{proof}

Our first application of the above result will be a local Lipschitz continuity estimate for the quantity $H^g( m , A, G)$.

\begin{proposition}\label{pr:10}
Denote
\begin{eqnarray*}
& &\Mjutska:=\max (\|m \|_{B^{-r}_{\infty\infty}},\|m^\prime  \|_{B^{-r}_{\infty\infty}}),\\
& &a :=\max (\| A\|_{B^t_{pp}\to B^{w}_{11}} ,\| A^\prime  \|_{B^t_{pp}\to B^{w}_{11}}),\\
& &v :=\max (\| G \|_{B^{\widetilde t}_{pp}\to B^{t}_{pp}} ,\| G^\prime \|_{B^{\widetilde t}_{pp}\to B^{t}_{pp}}).
\end{eqnarray*}
Then it holds that
\begin{eqnarray}\label{eq10'}
&&\| H^g( m, A, G)-H^g( m^\prime  , A^\prime , G^\prime )\|_{B^t_{pp}}\\
&\leq & (\| m-m^\prime \|_{B^{-r}_{\infty\infty}} +\| A-A^\prime\|_{B^{t}_{pp}\to B^{r}_{11}}+\| G-G^\prime
\|_{B^{\widetilde t }_{11}\to B^{t}_{11}})h(\Mjutska,a,v),\nonumber
\end{eqnarray}
where 
\ba
h(\Mjutska,a,v):=c\big(\Mjutska+a+(1+v)^q\big)\exp \left(
 a^{({2r/(w-r+2r/p^\prime)})}\Mjutska^{({2w/(w-r+2r/p^\prime)})} \right).\ea
Here $w$ and other indices are  as in Proposition \ref{pr1}.
\end{proposition}
\begin{proof}
Choose arbitrary $h^\prime \in B^{-t}_{p^\prime p^\prime}$ with $\| h^\prime \|_{B^{-t}_{p^\prime p^\prime}}=1.$ 
It is enough to estimate
the difference
\begin{eqnarray}\label{eq105'}
\expec|f(2,U)-f(1,U)|\nonumber,
\end{eqnarray}
where  for $x\in [1,2]$ and $U\in B^t_{pp}$ we have
\begin{eqnarray}\label{eq11'}
& &f(x,U)= \\ \nonumber
& &\bigg\langle h^\prime , g((G+(x-1)(G^\prime -G))U)
\exp \Big(-{1\over 2}\| (A+(x-1)(A^\prime -A))U\|^2_{L^2}+\nonumber\\
&&\phantom{hhhhh}+\langle  (A+(x-1)(A^\prime -A))U,m+(x-1)(m^\prime -m) \rangle_{}
\Big)\bigg\rangle\nonumber
\end{eqnarray}
It will be convenient to introduce for each $x\in [1,2]$ the notation $G_x=G+(x-1)(G^\prime -G),$ $A_x=A+(x-1)(A^\prime -A),$ and $m_x=m+(x-1)(m^\prime -m).$ As $g$ was assumed to be Lipschitz it follows that for each fixed
$U$ the function $x\mapsto f(x,U)$ is also Lipschitz.  Hence we may compute for almost every $x\in [1,2]$
\begin{eqnarray}\label{eq12'}
D_xf(x,U)
&=&\Big({d\over dx}\langle h^\prime , g(G_xU) \rangle +\langle h^\prime , g(G_xU) \rangle\big(-\langle (A^\prime -A)U,A_xU\rangle_{L^2}+\nonumber\\
&&\phantom{kk}+\langle A_xU,m^\prime -m\rangle_{}+\langle (A^\prime -A)U,m_x\rangle_{}\big)\Big)\cdot\nonumber\\
&&\phantom{kk}\cdot\exp \Big(-{1\over 2}\| A_xU\|^2_{L^2}+ \langle  A_xU,m_x \rangle_{} \Big).\nonumber
\end{eqnarray}
In order to estimate the right hand side  we observe first that our assumption on $g$ yields for almost every $x\in [1,2]$
\begin{eqnarray}\label{eq13'}
|{d\over dx}(\langle h^\prime , g(G_xU) \rangle)|&\leq& \|(G^\prime -G)U\|_{B^t_{pp}}(1+v\| U\|_{B^{\widetilde t}_{pp}})^q\nonumber\\
&\leq & \|G^\prime -G\|_{B^{\widetilde t }_{pp}\to B^t_{pp}}\|U\|_{B^{\widetilde t }_{pp}}(1+v\| U\|_{B^{\widetilde t}_{pp}})^q\nonumber\\
&\leq & \|G^\prime -G\|_{B^{\widetilde t }_{pp}\to B^t_{pp}}(1+v)^{q}(1+\| U\|_{B^{\widetilde t }_{pp}})^{q+1}.\nonumber
\end{eqnarray}
Moreover, by using the observation that $ B^{r}_{11}(\torus^d)\subset L^2(\torus^d)$ it follows that
\begin{eqnarray}\label{eq14'}
&&|-\langle (A^\prime -A)U,A_xU\rangle_{L^2}+\langle A_xU,m^\prime -m\rangle_{}+\langle (A^\prime -A)U,m_x\rangle_{}|\nonumber\\
&\leq & \| A_x\|_{B^t_{pp}\to L^2}\| A^\prime -A\|_{B^t_{pp}\to L^2}\| U\|^2_{B^t_{pp}} +\| m^\prime - m\|_{B^{-r}_{\infty\infty}}\| A_x\|_{B^t_{pp}\to B^{r}_{11}}\| U\|_{B^t_{pp}}
\nonumber\\
&&\phantom{mmmmmmmmm}+\| m_x\|_{B^{-r}_{\infty\infty}}\| A^\prime-A \|_{B^t_{pp}\to B^{r}_{11}}\| U\|_{B^t_{pp}}
\nonumber\\
&\leq &c(\Mjutska+a) (\| m^\prime -m\|_{B^{-r}_{\infty\infty}} +\| A^\prime-A \|_{B^{t}_{p ,p}\to B^{r}_{11}})(1+\|
U\|_{B^t_{pp}})^2. \nonumber
\end{eqnarray}
By combining the previous bounds and the obvious bound for $g(G_xU)$ we obtain
\begin{eqnarray}\label{eq15'}
|D_xf(x,U)| &\leq & c\big(\Mjutska+a+(1+v)^q\big)\big(\| m^\prime -m\|_{B^{-r}_{\infty\infty}}
+\\ \nonumber
&&+\| A^\prime-A \|_{B^{t}_{p ,p}\to B^{r}_{11}}+\| G^\prime -G\|_{B^{\widetilde t }_{pp}\to B^t_{pp}}\big) \cdot S_{q+2}(m_x , A_x).
\end{eqnarray}
The Fubini theorem allows us to compute
\begin{eqnarray}\label{eq18'}
\hspace{-3mm}&&\phantom{kuk}\expec |f(2,U)-f(1,U)| \leq  \int^2_1\Big( \expec |{d\over dx}(f(x,U))|\Big)\, dx \\
\hspace{-3mm}&\leq&\hspace{-1mm} c\big(\Mjutska+a+(1+v)^q\big)\big(\| m^\prime -m\|_{B^{-r}_{\infty\infty}}
+\| A^\prime-A \|_{B^{t}_{p ,p}\to B^{r}_{11}}+\| G^\prime -G\|_{B^{\widetilde t }_{pp}\to B^t_{pp}}\big)\nonumber\\
\hspace{-3mm}&&\phantom{kukkuu}\cdot \int^2_1 \expec S_{q+2}( m_x , A_x)\, dx.\nonumber
\end{eqnarray}
According to Proposition \ref{pr1} there is the uniform bound
\begin{eqnarray}\label{eq16'}
\expec S_{q+2}( m_x , A_x)\leq  c \exp \left( a^{({2r/(w-r+2r/p^\prime)})}
\Mjutska^{({2w/(w-r+2r/p^\prime)})} \right),
\end{eqnarray}
and the \refeq{eq10'} follows immediately by combining this with  \refeq{eq18'}.
 \hfill$\Box$
\medskip
\end{proof}

The local Lipschitz constant of a given map $r:E\to F$ between the Banach spaces $E$ and $F$ is defined as
\begin{eqnarray}\label{eq20'}
\mbox{Lip}(r)(x):=\limsup_{y\to x}\frac{\| r(y)-r(x)\|_F}{\| y-x\|_E}\quad \mbox{for}\;\; x\in E.\nonumber
\end{eqnarray}
In terms of this quantity the previous proposition states that
\begin{eqnarray*}
&&\mbox{Lip}(H^g)( m, A, G)\leq \\
&& \hspace{-10mm}c (1+\| G \|_{B^{\widetilde t}_{pp}\to B^{t}_{pp}})^q \exp \left(
 c(\| A  \|_{B^t_{pp}\to B^{w}_{11}})^{2r/(w-r+2r/p^\prime)}(\|m
\|_{B^{-r}_{\infty\infty}})^{{2w/(w-r+2r/p^\prime)}}\right) ,\nonumber
\end{eqnarray*}
where the linear factors containing   norms of $A$ and $m$ were absorbed in the exponential term. Above we consider $H^g$  as the
map
$$
H^g\, :\,\Big( B^{-r}_{\infty\infty}\oplus L(B^t_{pp}, B^{w}_{11})\oplus L(B^{\widetilde t}_{pp},B^{t}_{pp})\Big)\longrightarrow B^t_{pp}.
$$

Before we are able to give good estimates for the Lipschitz constant of the ratio $H^g/H^1$ we need a couple of auxiliary results.
\begin{lemma}\label{le:1'} Let $a>0$ and assume that $f$ and $F$ are non-negative  random variables with $\expec F <\infty $.
Then there is a constant $c$ depending only on $a$ so that
\begin{equation}\label{effkaava}
  \frac{\expec (fF)}{\expec F}\leq c\left(c+ \log (\expec \exp (f^a/2))+\log (\expec F^2)-2\log(\expec F )  \right)^{1/a}.
\end{equation}
\end{lemma}
\begin{proof} Let us consider the change of probability measure $\prob_F$ that is simply obtained by using $F$
as a weight and normalizing. Then the left hand side of (\ref{effkaava}) equals $\expec_Ff.$ By invoking the convex function
$\phi_a(x):= \exp((1+|x|)^a/a)$ and applying the Jensen inequality we obtain
$$
\phi_a(\expec_Ff)\leq \expec_F\phi_a(f)\leq \frac{(\expec (\phi_a(f))^2)^{1/2}(\expec F^2)^{1/2}}{\expec F},
$$
where the last inequality followed from the Cauchy-Schwartz inequality. By using the inequality $\frac{2}{a}(1+|x|)^a\leq c+c|x|^a$ and
applying the inverse $\phi_a^{-1}(x)=(a\log (x))^{1/a}-1$ a direct computation yields the inequality
\begin{eqnarray}\label{eq28'}
{\expec (fF)}/{\expec F} \leq  c\left(c+ \log (\expec \exp (cf^a))+\log (\expec F^2)-2\log(\expec F )  \right)^{1/a}.
\end{eqnarray}
Then (\ref{effkaava}) follows by applying \refeq{eq28'} on $(2c)^{-1/a}f$ in place of $f.$
\end{proof}
\begin{lemma}\label{ala} For any $t<s-n/p$ there is a constant
 $C=C(t)>0$  such that
$$
  H^1(m,A)\geq C(t)\exp(-{1\over 2}\| A\|_{B^{t}_{pp}\to L^2}^2),
$$
independently of $m.$
\end{lemma}
\begin{proof} Since the prior measure is invariant under the change of variables $U\to -U$
we obtain that
\begin{eqnarray}\label{eq31'}
&&\expec H^1( m , A)\nonumber\\
&=&{1\over 2}\expec\Big(\exp \big(-{1\over 2}\| AU\|^2_{L^2}+\langle  AU,m\rangle_{} \big)
+\exp \big(-{1\over 2}\| AU\|^2_{L^2}-\langle  AU,m\rangle_{} \big)\Big)\nonumber\\
&\geq &\expec\exp \big(-{1\over 2}\| AU\|^2_{L^2}\big).\nonumber
\end{eqnarray}
We may clearly take $C(t):=\prob (\{\| U\|_{B^t_{pp}}\leq 1\}) .$
 \hfill$\Box$
\medskip
\end{proof}

\begin{lemma}\label{le:2'} Let $\lambda >0.$  Under assumption $A$ and with $w$ as in Proposition \ref{pr1} there is
the estimate
$$
\frac{\expec S_\lambda( m , A)}{\expec S_0( m , A)}\leq c(1+\|m \|_{B^{-r}_{\infty\infty}})^{
\beta_1}(1+\| A \|_{B^t_{pp}\to B^{w}_{11}})^{\beta_2}
$$
where $\beta_1=\frac{2\lambda
w}{p(w-r+2r/p^\prime)}$ and $\beta_2=\max (\frac{2\lambda}{p}, \frac{2r\lambda}{p(w-r+2r/p^\prime)})$.
\end{lemma}
\begin{proof}
We apply Lemma \ref{le:1'} with the choice $a=p/\lambda,$ $f=(\frac{1}{2}(1+\|U\|_{B^{\widetilde t}_{pp}}))^\lambda$ and $F=\exp
\big(-{1\over 2}\| AU\|^2_{L^2}+\langle  AU,m\rangle_{} \big).$ One inserts
 the estimate for the quantity $\expec F^2 =I'$,  see \refeq{iprime},
 obtained in  (\ref{EQ: A1}) and  (\ref{EQ: A2})
in the proof of
Proposition \ref{pr1}. Moreover, a lower bound for $\expec F$ is given by  Lemma \ref{ala} and we recall that $ \expec
\exp (f^a/2))<\infty $  thanks to Lemma \ref{le1}. \hfill$\Box$
\medskip
\end{proof}

We are ready for the main result of this section.
\begin{theorem}\label{stability} Consider the ratio $H^g( m , A, G)/H^1( m , A, G)$ as a map
$$
  \frac{H^g}{H^1}\, : \, \Big( B^{-r}_{\infty\infty}\oplus L(B^t_{pp}, B^{w}_{11})\oplus L(B^{\widetilde t}_{pp},B^{t}_{pp})\Big)\longrightarrow B^t_{pp}.
$$
Under assumption $A$, and with $w$ as in Proposition \ref{pr1}, the local Lipschitz constant of this function satisfies
\begin{eqnarray}\label{eq40'}\\ \nonumber
&&\mbox{\rm Lip}(\frac{H^g}{H^1})( m , A, G)\leq c(1+\| G \|_{B^{\tilde t}_{pp}\to B^{t}_{pp}})^q(1+\|m
\|_{B^{-r}_{\infty\infty}})^{\gamma}(1+\| A \|_{B^t_{pp}\to B^{w}_{11}})^{\alpha},
\end{eqnarray}
where $\alpha =1+ (\frac{2q+8}{p})\max (1,(\frac{r}{w-r+2r/p^\prime}))$ and $\gamma= 1+(\frac{2q+8}{p})(\frac{w}{w-r+2r/p^\prime})$.
\end{theorem}
\begin{proof} Denote $\Mjutska=\|m \|_{B^{-r}_{\infty\infty}},$ $a=\| A \|_{B^t_{pp}\to B^{w}_{11}}$
and $v=\| G \|_{B^{\widetilde t}_{pp}\to B^{t}_{pp}}.$ Proposition \ref{pr:10} verifies that both $H^g$ and $H^1$ are locally Lipschitz. As a special case
we obtain that $\expec S_{q+2}(m,A)$
 is continuous with respect to variables $A$ and $m.$ Hence inequality \refeq{eq18'} yields the estimate
$$
\mbox{Lip}(H^g)( m, A, G)\leq c\big(\Mjutska+a+(1+v)^q\big) \expec S_{q+2}( m , A).
$$
The simple inequality (here $x_1, x_2$ are vectors and $y_1,y_2>0 $ are scalars)
$$
\| {x_1\over y_1}-{x_2 \over y_2 }\| \leq \frac{1}{y_2}\| x_2-x_1\|+\frac{|y_2-y_1 |}{y_1y_2}\,\| x_1 \|
$$
verifies that point-wise
\begin{eqnarray}\label{eq22'}
\mbox{Lip}(\frac{H^g}{H^1})\leq \frac{\mbox{Lip}(H^g)}{H^1}+\frac{\|H^g\|_{B^{t}_{pp}}\mbox{Lip}(H^1)}{(H^1)^2}.\nonumber
\end{eqnarray}
Observe that by assumption A we have  $\|H^g\|_{B^{t}_{pp}}\leq c\expec S_q.$ By combining the previous estimates and remembering that
$H^1=\expec  S_0$ we thus obtain
$$
\mbox{Lip}(\frac{H^g}{H^1})\leq c\big(\Mjutska+a+(1+v)^q\big) \Big(\frac{\expec S_{q+2}}{\expec S_0}\Big)+ c\Big(\frac{\expec S_{q}}{\expec
S_0}\Big) c\big(\Mjutska+a\big)\Big(\frac{\expec S_{2}}{\expec S_0}\Big).
$$
The statement then follows by three applications of Lemma \ref{le:2'}. \hfill$\Box$
\medskip
\end{proof}

\begin{remark}{\rm The main content of Theorem \ref{stability} is a stability estimate that grows polynomially in the norm of the measurement $m.$ We have not striven for the optimal exponents in the above computations.
Moreover, one should observe that in the case $p>1$ it is possible to choose $w=r$ in Theorem \ref{stability}, which corresponds to the
weakest condition on smoothness of $A$. If $p\in (1,2)$ this yields the exponents $\alpha=\gamma=1+p^\prime(q+4)/p$. In the important
special case of $p=1$ one is forced to choose $w>r.$ On the other hand, the choices $w>r$ yield considerably smaller exponents for all
small values $p\geq 1$: in the limit $w\to\infty$ one has $\alpha\to 1+ (\frac{2q+8}{p})$, and $\gamma\to 1+ (\frac{2q+8}{p})$. }

\end{remark}

\section{Convergence results for Besov priors}\label{sec:BesovConv}

In the present section we assume that indices $t,\widetilde t ,r $ and the quantities $m$, $A$, and $g$ satisfy Assumption A. We
take $T_n:B^{\widetilde t }_{pp}(\torus^d)\to B^{\widetilde t }_{pp}(\torus^d)$, $n\geq 1$, defined by the familiar truncation
\begin{equation}\label{truncFromula}
  T_n\left(\sum_{\ell=1}^\infty c_\ell \psi_\ell\right) =  \sum_{\ell=1}^n c_\ell \psi_\ell
\end{equation}
discussed above and in Appendix \ref{app:wavelet}. Then $T_n\to I$ strongly in $L(B^{\widetilde t }_{pp})$ as $n\ra \infty$. We consider proper linear discretizations $U_n=T_nU$ of the random variable $U$.

By standard compact imbedding results it follows that $(I-T_n)\to 0$ as $n\ra \infty$ in the operator norm topology of $L(B^{\widetilde t }_{pp},B^{t}_{pp})$, and $(I-P_k)\to 0$ as $k\ra \infty$ in the operator norm topology of $L(B^{\widetilde r }_{11},B^{r}_{11})$ with $\widetilde r >r$,  see Appendix \ref{app:wavelet}. {\mntext Next we formulate the assumptions in quantitative
terms.}
\medskip

\noindent {\bf Assumption B:}{\it \;\; Let $\widetilde r>r$ and assume that $A:B_{pp}^{t}(\torus^d)\to B_{11}^{\widetilde r }(\torus^d)$.
Define $T_n$ by (\ref{truncFromula}); then we have
\begin{equation}\label{assu 3c}
  \|I-T_n\|_{L(B^{\widetilde t }_{pp}, B^{t}_{pp})}<\eta_2(n) \qquad \mbox{ for }n\geq 1
\end{equation}
with $\lim_{n\to \infty}\eta_2(n) =0$. Further, let $P_k:B^{\widetilde r }_{pp}(\torus^d)\to B^{\widetilde r }_{pp}(\torus^d)$ be bounded linear projections
with $k$-dimensional range satisfying
\begin{equation}\label{assu 3}
  \|I-P_k\|_{L(B^{\widetilde r }_{11}, B^{r}_{11})}<\eta_1(k), \qquad \mbox{ for }k\geq 1
\end{equation}
with  $\lim_{k\to \infty}\eta_1(k) =0$. Assume that $\|T_n\|_{L(B^{\widetilde t }_{pp})}\leq C$ and $\|P_k\|_{L(B^{\widetilde r
}_{pp})}\leq C$ for all $k, n\geq 1$. The measurements are assumed to be uniformly bounded by a constant $\Mjutska$:
\begin{eqnarray}
\label{assu 8} \|m\|_{B^{-r}_{\infty\infty}}\leq \Mjutska\quad \mbox{and}\quad
 \|m_k\|_{B^{-r}_{\infty\infty}}\leq \Mjutska\quad \mbox{for all} \;\;k\geq 1.
\end{eqnarray}
Finally, it is  assumed  that
\begin{eqnarray}
 \label{assu 10} \eta_3(k):=\|m_{k}-m\|_{B^{-r}_{\infty\infty}}\hbox{ satisfies }
\lim_{k\to \infty}\eta_3(k)=0.
\end{eqnarray}
}

Let $\Meas_{kn}=P_kAU_n+\Noise$. As before we {define the reconstructor
$\mathcal{R}_{\Meas_{kn}}(g(U_n)|m_k)$ corresponding to model $\Meas_{kn}$ at $m_k$} as

\begin{eqnarray}\label{deterministic, spesific2}
  \mathcal{R}_{\Meas_{kn}}(g(U_n)|m_k)= \frac{H_{{U_n,\Meas_{kn}}}^g(m_k)}{H^1_{{U_n,\Meas_{kn}}}(m_k)}.
\end{eqnarray}
The following result yields a quantitative convergence result for Besov priors:
\begin{theorem}\label{th:besovconvergence} Let Assumptions A and B hold.
Denote $\gamma =1+(\frac{2q+8}{p})(\frac{w}{w-r+2r/p^\prime})$. There is a constant $C'=C'(A,g,t,r,w,p)$  such that
\begin{eqnarray*}
\|\mathcal{R}_{\Meas_{kn}}(U_n|m_k)-\mathcal{R}_M(U|m)\|_{B^t_{pp}} \leq C^\prime (1+K)^{\gamma} [\eta_1(k)+\eta_2(n)+\eta_3(k)].
\end{eqnarray*}
Moreover,  the limits
$$
  \lim_{n\to \infty}\mathcal{R}_{\Meas_{kn}}(g(U_n)| m_k ) \quad\mbox{and}
  \quad   \lim_{k\to \infty}\mathcal{R}_{\Meas_{kn}}(g(U_n)| m_k )
$$
exist for a fixed value of $k$ (resp. $n$).
\end{theorem}
\begin{proof}
The statement is an immediate consequence of the estimates obtained in the previous section: observe that in our notation
$$
  \frac{H_{{U_n,\Meas_{kn}}}^g(m_k)}{H^1_{{U_n,\Meas_{kn}}}(m_k)} =\frac{H^g(m_{k}, P_kAT_n,
  T_n)}{H^1(m_k, P_kAT_n, T_n)}.
$$
Moreover, apart from a possible change of the uninteresting constants, the Lipschitz bound given by Theorem \ref{stability} on the quantity
$H^g/H^1$ is uniform on any bounded neighbourhood of $m,$ $A$ and {\mntext $G$,
where $G$ is the identical embedding} $\mbox{Id}: B^{\widetilde t}_{pp}\to B^{t}_{pp}$. Note
also that by our assumptions  $(A-P_kAT_n)\to 0$ in the operator norm topology of $L(B^{t}_{pp},B^{r}_{11})$.  The first claim
now follows from definitions and Theorem \ref{stability}. The last statements follow immediately
by the same reasoning.
 \hfill$\Box$
\medskip
\end{proof}

We emphasize  that Theorem \ref{th:besovconvergence} covers the highly interesting value $p=1$ despite the failure of assumption
(\ref{abc}) of the general theory in that case.


Let us revisit the deblurring example of Section \ref{sec:Besov_deblur}, where $p=1$, $s=1$, $d=2$ and $A$ has a smooth kernel. Often it is possible to take the
projection $P_k$ related to the measurement device to be truncation of the wavelet expansion to the {\ntext first $k$ terms} analogously to
(\ref{truncFromula}). (For instance, $P_k$ might measure local averages at a grid of points and then compute discrete wavelet transform by
convolutions with finite filters.)

Then Theorem \ref{th:besovconvergence} yields the convergence of reconstructors
and a result analogous to Theorem \ref{Corollary-teoreema} with an explicit convergence
speed.

\begin{corollary}\label{co:10}  Let Assumptions A and B hold. Let $p=s=1$ and assume that
$A:{\mathcal D}^\prime (\torus^2)\to C^\infty(\torus^2)$ is a bounded linear operator. Let $P_k$ and $T_n$ be truncations to $k$ and $n$ first terms in wavelet expansion, respectively. Moreover, let $t<\tilde t<-1$, $r>r_1>1$, $\lambda>11$, and $\tau>0$. Then
 \medskip

\noindent (i) There is a constant $c_0=c_0(A,r,t,\lambda)>0$ such that the reconstructors satisfy
\begin{eqnarray*}
  \|\mathcal{R}_{\Meas_{kn}}(U_n| m_k)
-\mathcal{R}_M(U|m)\|_{B^{t}_{11}(\torus^2)}
 &\leq&
 c_0 (1+K)^\lambda [\eta_1(k)+\eta_2(n)+\eta_3(k)].
\end{eqnarray*}

\medskip

\noindent (ii) Let $u=U(\omega_0)$,  $\e=\Noise(\omega_0)$,
$\omega_0\in \Omega$ be realizations of the random variables $U$ and $\Noise$, and
\ba
m = A u + \e,\quad  m_k = A_k u+P_k\e,
\ea
be the realizations of the measurements (\ref{basicequation A1}) and (\ref{eq: measmodel_k k}), respectively.
Moreover, assume that $m \in B^{-r_1}_{11}(\torus^2)$. Then
there is $C>0$ independent of $n$ and $k$ so that
the reconstructors  defined in Theorem \ref{g-theorem}
for  measurements
 (\ref{eq: measmodel_k k}) and (\ref{basicequation A1}) satisfy
\begin{eqnarray}
  \|\mathcal{R}_{M_{kn}}(U_n|m_k)
-\mathcal{R}_M(U|m)\|_{B^{t}_{11}(\torus^2)}
 &\leq&
C[k^{-\tau}+n^{-(\widetilde t- t)/2}].
\end{eqnarray}

\end{corollary}
\medskip

Note that in (ii) we have  $m \in B^{-r_1}_{11}(\torus^2)$ for $\prob$-a.e.\ $\omega$.
\medskip

\begin{proof} Claim (i) follows directly from Theorem \ref{th:besovconvergence} when we use $q=1$ and so large $w$ that $\lambda\geq \gamma$.
As $A$ is an infinitely smoothing operator, we can choose above any $-\infty<t<\widetilde t<-1$ and $\tilde r> r+2\tau>r_1+4\tau$. Moreover, since
$m \in B^{-r_1}_{pp}(\torus^2)$, we can take  above $\eta_1(k)=c_1 k^{-(\widetilde r-r)/2}$, $\eta_2(n)=c_2 n^{-(\widetilde t- t)/2}$, and
$\eta_3(k)=c_3 k^{-(r-r_1)/2}$ where $c_1,c_2$ and $c_3$ depend on $\omega_0$
and the parameters $\widetilde r,r,r_1,\widetilde t,t,$ but not on $k$ or $n$. As the projections $P_k$ satisfy the conditions
(\ref{Pk 1}) and (\ref{Pk 2}),
the assertion (ii) follows from (i) and Lemma \ref{reconstructors coincide}.
 \hfill$\Box$
\end{proof}
\medskip

\appendix

\section{Besov spaces and wavelets}\label{app:wavelet}

Let $\tilde{\psi}$ and $\tilde{\phi}$ be compactly supported wavelet and scaling function suitable for multi-resolution analysis  of
smoothness $C^r$ in $L^2(\R)$.

Following Daubechies \cite[section 9.3]{Daubechies} we construct a wavelet representation for periodic functions in $\R$ with period $1$;
in other words, for functions on the one-dimensional torus $\T^1$. Set
\begin{eqnarray*}
  \phi_{j,k}(x) &=& \sum_{\ell\in\Z}\tilde{\phi}(2^j(x+\ell)-k),\\
  \psi_{j,k}(x) &=& \sum_{\ell\in\Z}\tilde{\psi}(2^j(x+\ell)-k).
\end{eqnarray*}
We use in the following the subspaces of $ L^2(\T^1)$, 
$$
  V_j:=\overline{\mbox{span}\{\phi_{j,k}\,|\, k\in\Z\}},\qquad
  W_j:=\overline{\mbox{span}\{\psi_{j,k}\,|\, k\in\Z\}}.
$$
It turns out that $V_j$ are spaces of constant functions for $j\leq 0$. Thus we have a ladder $V_0\subset V_1 \subset V_2\subset \cdots$ of
multiresolution spaces satisfying
\begin{equation}\nonumber
  \overline{\bigcup_{j\geq 0} V_j} = L^2(\T^1).
\end{equation}
Further, we denote the successive orthogonal complements of $V_j$ in $V_{j+1}$ by $W_j$ for $j\geq 0$. Then we have orthonormal bases
\begin{eqnarray*}
  &&\{\phi_{j,k} \,|\, k=0,\dots,2^j-1\}\mbox{ in } V_j,\\
  &&\{\psi_{j,k} \,|\, k=0,\dots,2^j-1\}\mbox{ in } W_j.
\end{eqnarray*}

Following Meyer \cite[section 3.9]{meyer} we define a wavelet basis for periodic functions in $\R^d$; in other words, for functions on the
torus $\T^d$. Let $E$ denote the set of $2^d-1$ sequences $\nu=(\nu_1,\nu_2,\dots,\nu_d)$ of zeroes and ones, excluding the sequence
$(0,0,\dots,0)$. Define for $\nu\in E$ and $j\geq 0$ the wavelets
$$
  \psi_{j,k}^\nu(x):= 2^{dj/2}\psi^{\nu_1}(2^jx_1-k_1)\dots\psi^{\nu_d}(2^jx_d-k_d)
$$
with the convention that $\psi^0=\phi$ and $\psi^1=\psi$, and integer-valued components of vector $k$ ranging over
$$
  0\leq k_i \leq 2^j-1 \quad \mbox{ for all }\quad i=1,2,\dots,d.
$$

The functions $\psi_{j,k}^\nu(x)$ constitute an orthonormal basis for $L^2(\T^d)$. Let us renumber the above basis functions using just one
integer $\ell=1,2,\dots$. First, $\ell=1$ corresponds to the scaling function $\phi(x_1)\dots\phi(x_n)$. The remaining numbering is done
scale by scale; that is, we first number wavelets with $j=0$, then wavelets with $j=1$, and so on. The $2^d-1$ indices $\nu\in E$ are
naturally numbered by thinking them as binary representation of integers. The exact ordering of all $2^{jd}$ translations corresponding to
a fixed $j$ can be chosen arbitrarily. This leads to a numbering of the following type:
\begin{eqnarray*}
  \mbox{scale }j=0: && \ell=2,\dots,2^d,\\
  \mbox{scale }j=1: && \ell=2^d+1,\dots,2^{2d},\\
  \mbox{scale }j=2: && \ell=2^{2d}+1,\dots,2^{3d},\\
  \vdots &&\vdots
\end{eqnarray*}

According to Meyer \cite[Section 6.10]{meyer}, we can characterize periodic Besov space functions using these wavelets. Namely, the series
$$
  f(x) = \sum_{\ell=1}^\infty c_\ell \psi_\ell(x)
$$
belongs to $B^{s}_{pq}(\T^d)$ if and only if
$$
  2^{js}2^{dj(\frac{1}{2}-\frac{1}{p})}\left(\sum_{\ell=2^{jd}}^{2^{(j+1)d}-1}|c_\ell|^p \right)^{1/p}\in\ell^q(\N)
$$
We always assume that $r$ is large enough for providing bases for Besov spaces with smoothness $s$. The case $q=p$ is especially relevant to us, and we obtain  the equivalent
norm \beqla{normi} \|\sum_{\ell=1}^\infty c_\ell \psi_\ell(x) \|_{B^{s}_{pp}(\T^d)} :=\left(\sum_{\ell=1}^\infty
\ell^{(ps/d+p/2-1)}|c_\ell|^p \right)^{1/p}. \eeq We use the above quantity as the definition of the Besov norm $\| \cdot
\|_{B^{s}_{pp}(\T^d)}$ for generalized functions  on $\torus^d$, $s\in\real$ and $p\in [1,\infty ].$ In case $p=\infty$ the
definition must be understood as follows:
$$
\|\sum_{\ell=1}^\infty c_\ell \psi_\ell(x) \|_{B^{s}_{pp}(\T^d)} :=\sup_{\ell\geq 1} \ell^{(s/d+1/2)}|c_\ell| .
$$

It follows that ${B^{s}_{pp}(\T^d)}$ is isometrically isomorphic to the sequence space $\ell^p,$ and by the simplicity of the norm it is
easy to control the basic properties of the spaces. Especially,
there is an  embedding (an easy corollary of the H\"older inequality),
 \beqla{besovembedding} B^{s_1}_{p_1p_1}(\T^d)\subset
{B^{s_2}_{p_2p_2}(\T^d)}\quad \mbox{if and only if} \quad s_1-d/p_1\geq s_2-d/p_2,
\eeq and it is easy to verify that this embedding is compact if $s_1-d/p_1>
s_2-d/p_2.$ The dual of $B^{s}_{pp}(\T^d)$ is $B^{-s}_{p^\prime p^\prime}(\T^d)$ if $p\in (1,\infty )$, and $p^\prime$ stands for the dual
index: \; $1/p+1/p^\prime=1.$ The duality is with respect to the standard duality bracket: if $f=\sum_{\ell=1}^\infty f_\ell \psi_\ell$ and
$g=\sum_{\ell=1}^\infty g_\ell \psi_\ell$ are finite sums, then
$$
\langle f,g\rangle =\int_{\torus^d}f(x)g(x)\, dx=\sum_{\ell= 1}^\infty f_\ell g_\ell.
$$
We finally observe that natural  bounded linear projection operators on $B^{s}_{pp}(\T^d)$ are
obtained by setting
$$
  T_n f (x) = \sum_{\ell=1}^n c_\ell \psi_\ell(x).
$$
These projections work at the same time for all the spaces and are contractions.
Moreover, if $p<\infty ,$ then $\lim_{n\to\infty} \| f-T_nf\|_{B^{s}_{pp}(\T^d)}=0$
for all $f\in {B^{s}_{pp}(\T^d)}.$

\section{{\ntext Examples of limits of finite-dimensional random variables}}\label{sec:discretizationexamples}
{\ntext
Here we illustrate difficulties related to finite-dimensional models and their possible convergence to an infinite-dimensional continuum model.
Unless otherwise stated, we work on a circle $\T^1$, or equivalently, the interval $[0,1]$ with end points identified. Let $u\in L^2(\T^1)$. We consider two measurements
\ba
A^{(1)}u=\int_{\T^1} u(t)\,dt,\quad A^{(2)}u=u(\frac 12)
\ea
and the corresponding measurement models
\beq\label{measurement model type 1}
& &M^{(1)}_n= A^{(1)}U_n+\Noise,\quad M^{(1)}= A^{(1)}U+\Noise,\\
& &M^{(2)}_n= A^{(2)}U_n+\Noise,\quad M^{(2)}= A^{(2)}U+\Noise \label{measurement model type 2}
\eeq
where $U_n$ is a finite-dimensional random variable, $U$ is a random variable in an infinite-dimensional function space, and $\Noise$ is a normalized Gaussian random variable independent of $U$ and $U_n$. We consider various examples of $U_n$ and study whether some random variable $U$ could be considered as a limit of $U_n$ as $n\to \infty$, and whether models (\ref{measurement model type 1}) or (\ref{measurement model type 2}) make sense in the limit.

In all examples below, $X_j^n$, $j,n\in \N$ are independent real-valued normalized Gaussian random variables.
\medskip

{\bf Example 1.} (``Non-proper discretization of white noise'') Let $I(n,j)=(\frac {j-1}n,\frac {j}n]$, $j=1,\dots,n$ and
$\chi^n_j(t)=\chi_{I(n,j)}(t)$ be the indicator functions of intervals $I(n,j)$.
Let
\beq\label{new example 1}
U_n(t)=\sum_{j=1}^{n} a^n_jX_j^n\chi^n_j(t), \quad t\in \T^1
\eeq
where $a^n_j>0$ are parameters.

For a fixed $n$, with an {\it ad hoc} choice $a^n_j=1$, the functions $U_n$ could be considered as an interesting {\nnntext random signal}. However, for any function $\phi\in L^2({\T^1})$
\ba
\lim_{n\to \infty}\int_{\T^1}U_n(t)\phi(t)\,dt=0\quad\hbox{in distribution},
\ea
and thus $U_n$, considered as $L^2({\T^1})$ valued random variables, converge to zero {\nnntext weakly} in distribution as $n\to \infty$. Taking $U=0$  we see that the measurements $M^{(1)}_n$ converge in distribution to $M^{(1)}$. Concerning measurement (\ref{measurement model type 2}) {\nnntext we notice the $U_n(\frac 12)\sim N(0,1)$, but $U(\frac 12)\sim N(0,0)$.} Thus $M^{(2)}_n\sim N(0,2)$ for all $n$, but for $U=0$ we have $M^{(2)}\sim N(0,1)$. This shows that measurement models (\ref{measurement model type 2}) do not behave nicely with the choice $a^n_j=1$. In this example $U_n$ are not proper linear discretizations of $U$ in any Banach space, since the second condition in Definition \ref{def:linear_discr} is violated.

\medskip

{\bf Example 2.}  (``Proper discretization of white noise'')
Consider random variables $U_n$ defined in (\ref{new example 1})
with parameters $a^n_j=n^{1/2}$. This choice is motivated by the fact that the functions $n^{1/2}\chi_{I(n,j)}$
are orthonormal in $L^2(\T^1)$.
Then, for $\phi\in C^\infty({\T^1})$
\begin{equation}\label{esimeqqq}
\lim_{n\to \infty} \int_{\T^1}U_n(t)\phi(t)\,dt=\bra U,\phi\cet \quad\hbox{in distribution},
\end{equation}
where $U$ is Gaussian white noise in $L^2({\T^1})$. 
This  actually holds also when $\phi\in H^1 (\T^1).$ 
Let $Q_n$ be the $L^2(\T^1)$-orthogonal projection
on the subspace of functions that have constant value on the intervals $I(n,j)$.
Then $\bra U_n,\phi\cet$ appearing on the left hand side of (\ref{esimeqqq}) is a real-valued Gaussian random variable with covariance 
\ba
\sum_{j=1}^n |\langle \phi,n^{1/2}\chi_{I(n,j)}\cet|^2=\| Q_n\phi\|_2^2.
\ea 
As $n\to\infty$ this tends to the value $\|\phi\|_2^2$, as is easily seen by
first approximating $\phi$ by smooth functions. Hence $U_n\to U$ weakly in distribution in the
space $H^{-1}(\T^1).$

Moreover, we note that
$U_n$ has the same distribution as the random variable $T_nU,
$
where  $T_n: H^{-1}(\T^1)\to H^{-1}(\T^1)$ is the linear operator 
$T_nv=\sum_{j=1}^ n
\bra v,\phi^n_j\cet\, n^{1/2}\chi_{I(n,j)}$, 
where $(\phi^n_j)_{j=1}^ n$ is any orthonormal sequence in $L^2(\T^1)$ consisting
of elements of $H^1(\T^1)$.
We have thus verified that the variables $U_n$ are proper linear discretizations of $U$ in the space
$H^{-1}(\T^1)$, according to Definition \ref{def:linear_discr}.

Now the measurement $M^{(1)}$ is well defined and $M^{(1)}_n$ converge in distribution to $M^{(1)}$ as $n\to \infty$.
However, we have $M^{(2)}_n\sim N(0,n)$ and thus measurements $M^{(2)}$ do not converge in distribution
as $n\to \infty$. This is related  to the fact that the white noise
$U$ is a well defined $H^{-1}({\T^1})$ valued Gaussian random variable, and the constant
function 1 is in the dual of the space  $H^{-1}({\T^1})$ but the point evaluation $u\mapsto u(\frac 12)$
does not define a bounded linear operator in  $H^{-1}({\T^1})$.

\medskip

\begin{remark}\label{re:examples} {\rm In Example 2 above one may verify that it is possible to choose the $L^2(\T^1)$-orthonormal sequence $(\phi^n_j)_{j=1}^n$ so that the norm $\|T_n \|_{H^{-1}({\T^1})\to H^{-1}({\T^1})}$ remains uniformly bounded  for all $n$. However, this choice is somewhat complicated. A way to construct discretizations with this property (and such  that they fall in the scope of the basic theory developed in Section \ref{sec:general}) is to apply suitable  finite dimensional approximations of identity that are uniformly bounded simultaneously on both $H^{-1}({\T^1})$ and $L^2({\T^1}).$ E.g., one may  truncate Fourier series or apply basis projections corresponding to  a wavelet basis (set $p=2$ in Section \ref{sec:BesovConv}). Details of these comments will be considered elsewhere. Similar remarks apply to Example 3 below.}
\end{remark}

{\bf Example 3.} (``Discretization of the Gaussian smoothness prior'')  Choose continuous functions $\eta^n_j:{\T^1}\to \R$ so that they are affine on intervals $I(n,j)$, $j=1,\dots,n$ and that 
$(\eta^n_j)_{j=1}^n$ are orthonormal in $H^1({\T^1})$. Let
\beq\label{new example 3}
U_n(t)=\sum_{j=1}^{n} b^n_jX^n_j\eta^n_j(t),
\eeq
where $b^n_j>0$ are parameters. We choose $b^n_j=1$.
Then, for $\phi\in L^2({\T^1})$
\begin{equation}\label{eq:gash}
\lim_{n\to \infty} \int_{\T^1}U_n(t)\phi(t)\, dt=
 \int_{\T^1}U(t)\phi(t)\,dt\quad\hbox{in distribution},
\end{equation}
where $U$ is a Gaussian random variable 
in $L^2(\T^1)$ having zero expectation and covariance operator $(I-\Delta)^{-1}$,
that is, $U$ is the one-dimensional Gaussian
smoothness prior.
Let $Q_n$ be the $H^1 (\T^1)$-orthogonal projection onto the subspace $Y_n$ spanned by
$\eta^n_j$, $j=1,2,\dots,n$. 
Analogously to Example 2, one can see that $U_n$ have the same distribution as 
random variables $T_nU$ where
where  $T_n: L^2(\T^1)\to L^2(\T^1)$ is the linear operator 
$$
T_nv=\sum_{j=1}^ n
\bra v,(I-\Delta)\phi^n_j\cet\, \eta^n_j,
$$ 
where $(\phi^n_j)_{j=1}^ n$ is any orthonormal sequence 
in $H^1(\T^1)$ consisting
of elements of $H^2(\T^1)$. Thus 
 $U_n$ are proper linear discretizations of $U$ in $L^2(\T^1).$ 

Let us next consider $\phi\in H^{-1}(\T^1)$. 
Due to the formula (\ref{new example 3}), the $(H^{1}(\T^1)\times H^{-1}(\T^1))$-duality 
$\bra U_n(\omega),\phi\cet$
defines a real-valued Gaussian random variable with covariance
\beq\label{eq: Appendix B 1}
\expec (\bra U_n,\phi\cet^2)&=&
\expec (\sum_{j,k=1}^n b^n_jX_j\,b^n_kX_k\, \bra \eta^n_j,\phi\cet
\bra \eta^n_k,\phi\cet)
\\ \nonumber
&=&
\sum_{j=1}^n \bra \eta^n_j,\phi\cet^2
\\ \nonumber
&=&  
\sum_{j=1}^n \bra \eta^n_j,(I-\Delta)^{-1}\phi\cet_{H^1(\T^1)}^2\\
\nonumber
&=&\|Q_n(I-\Delta)^{-1}\phi\|_{H^{1}(\T^1)}^2. 
\eeq
The kernel of the covariance of operator of $U$ in $L^2(\T^1)$,
that is, the function 
\ba
G(t,t')=\expec (U(t)U(t')), \quad t,t'\in \T^1
\ea 
is Green's 
function of the differential operator $-\frac {d^2}{dt^2}+1$, 
\ba
(-\frac {d^2}{dt^2}+1)G(t,t')=\delta(t-t').
\ea
This implies that
$\expec (U(t)U(t'))$ is 
Lipschitz smooth on $\T^1\times \T^1$. As $U$ is Gaussian, one can see using 
e.g. \cite[Theorem 3.23]{kall}  that
 the values of the Gaussian smoothness prior $U$ are almost
surely in 
any H\"older space $C^\alpha(\T^1)$ with
$\alpha <1/2$. Thus, after fixing $\alpha\in (0,\frac 12)$, we can consider $U$ also
as a Gaussian $C^\alpha(\T^1)$-valued random variable satisfying
\beq\label{eq: Appendix B 2}
\expec (\bra U,\phi\cet^2)=\|(I-\Delta)^{-1}\phi\|_{H^{1}(\T^1)}^2
\eeq 
for every $\phi$ in the dual space of $C^\alpha(\T^1)$.
Note that as $H^1(T^1)\subset C^\alpha(\T^1) $, we have
$( C^\alpha(\T^1))'\subset H^{-1}(T^1)$.
Thus, since the projectors
$Q_n$ converge strongly to identity
in $H^{1}(\T^1)$ as $n\to \infty$, we can use  (\ref{eq: Appendix B 1}), 
(\ref{eq: Appendix B 2}), and the fact that the ranges of the operator $T_n$ used
above are in $H^1(\T^1)\subset C^\alpha(\T^1)$, and infer that
that $U_n$ are  proper linear discretizations of $U$ in $C^\alpha(\T^1)$, too. 
Because of this the measurements
$M^{(1)}$ and $M^{(2)}$ are well defined and one can see that the measurements
$M^{(1)}_n$ and $M^{(2)}_n$  converge in distribution to $M^{(1)}$ and $M^{(2)}$, respectively,
 as $n\to \infty$.

\medskip

{\bf Example 4.}  (``Discrete total variation priors'')  Let us next consider an example on interval $I=[0,1]$ (i.e., the end points
are not identified) and let  $\theta^n_j:[0,1]\to \R$ be continuous functions with are
affine functions on intervals $I(n,j)$, $j=1,\dots,n$, vanishing at $t=0$ and
$t=1$ such that
$\theta^n_j,$ $j=1,2,\dots, n-1$ are orthonormal in $H^1(I)$.
Let
\beq\label{new example 4}
U_n(t)=\sum_{j=1}^{n-1} Z_j^n\eta^n_j(t),
\eeq
where $Z^n=(Z^n_1,Z^n_2,\dots,Z^n_{n-1})$ is a $\R^n$ valued random variable
having the probability density function
\ba
\pi_{Z^n}(z_0,\dots,z_{n-1})=c_n \exp\left(-a_n\|\frac d{dt}(\sum_{j=1}^{n-1} z_j\eta^n_j(t))\|_{L^1(I)}\right)
\ea
where $a_n>0$ is a parameter and $c_n$ is a normalization constant of the probability
density function. The distribution of the random variables $U_n$ are sometimes
called the discrete total variation prior.
By \cite{LS}, the random variables $U_n$ converge in distribution if
$a_n=n^{1/2}$ but then the limit $U$ is a {\it Gaussian} random variable.
As a Gaussian distribution stays Gaussian in a linear transformation,
we see that
in this example $U_n$ are not proper linear discretizations of $U$ in
any Banach space.
\medskip

{\bf Example 5.}  (``Discretization of the two-dimensional Gaussian smoothness prior'')  Let us consider also higher dimensional example
analogous to Example 3. Let $L(n,j,k,1),L(n,j,k,2)\subset
\T^2$ be disjoint triangles so that their union is the square
$I(n,j)\times I(n,k)\subset \T^2$, $j,k=1,\dots,n$.
Let
 $\zeta^n_{j,k,l}:{\T^2}\to \R$, $j,k=1\dots,n$, $l=1,2$ be continuous functions on $\T^2$
which are
affine functions on triangles $L(n,j,k,l)$
such that
$\zeta^n_{jkl},$ $j,k=1,2,\dots, n$, $l=1,2$ are orthonormal in $H^1({\T^2})$. Let
\beq\label{new example 5}
U_n(t_1,t_2)=\sum_{j,k=1}^{n}\sum_{l=1}^2 b^n_{jkl}X^n_{jkl}\zeta^n_{jkl}(t_1,t_2),\quad (t_1,t_2)\in \T^2
\eeq
where $X^n_{jkl}\sim N(0,1)$ are independent and  $b^n_{jkl}>0$ are parameters.
Let $b^n_{jkl}=1$. Then the $U_n$
can be considered as Gaussian random variables in $H^{-1}(\T^1)$
that converge weakly in distribution to a $H^{-1}(\T^1)$ valued Gaussian random variable $U$
having zero
expectation  and the covariance operator $(I-\Delta)^{-2}$ in  $H^{-1}(\T^1)$.
This means that $U$ is the two-dimensional Gaussian smoothness prior.

Let now $u\in L^2(\T^1)$. We consider two measurements
\ba
A^{(1)}u=\int_{\T^2} u(t)\,dt,\quad A^{(2)}u=u(\frac 12,\frac 12)
\ea
and define models $M^{(1)}_n,M^{(1)}$ and $M^{(2)}_n,M^{(1)}$
as in (\ref{measurement model type 1}) where $\Noise$ is
Gaussian white noise in $L^2(\T^2)$.
Then, the measurement
$M^{(1)}$ is well defined and the measurements
$M^{(1)}_n$  converge in distribution to $M^{(1)}$.
However, one can show that $M^{(2)}_n\sim N(0,\sigma_n^2)$, where
$\sigma_n\to\infty$   as $n\to \infty$ and we see
that the measurements  $M^{(2)}_n$ do not converge in distribution
 as $n\to \infty$.
In this example one can verify that the $U_n$ are proper linear discretizations of $U$ in
$H^{-1}(\T^2).$

The fact the point value measurements $M^{(2)}_n$ do not converge is
related to the fact that the two-dimensional Gaussian smoothness prior
has, formally speaking, the covariance function
$\expec(U(t)U(s))=G(t,s)$, $t,s\in \T^2$ is the Green's function of the operator $-\Delta+I$ and has thus on the diagonal $t=s$ a logarithmic
singularity. This fact is extensively used in quantum field theory
in the study of the free Gaussian field, a random field that is very
similar to the Gaussian smoothness prior \cite{Shef}.
}

\section{On the domain of reconstructors}\label{sec:recDomain}

Considering formula (\ref{fin dim formula}) in the infinite-dimensional case, we meet the difficulty that realizations of $M$ belong to $Z$
only with probability zero. Therefore the function  $m\mapsto \mathcal{R}_M(U| m)$ should be defined in some larger set than  $Z$.

A generalized definition of a reconstructor is as follows:
\begin{definition}\label{def: conditional mean function2}
The deterministic function $\mathcal{R}_M(U|\,\cdotp\,):S_0\to Y$, where $m\mapsto \mathcal{R}_M(U|\, m)$, defined in a Borel-measurable
subspace $S_0\subset S$ is {\em reconstructor} of $U$ with measurement $M$ if $M(\omega)\in S_0$ almost surely and
\begin{eqnarray}\label{det. CM estimate}
  \expec(U|\M)(\omega)=  \mathcal{R}_M(U|M(\omega))\quad \hbox{almost surely.}
\end{eqnarray}
The quantity $E(g(U)|\,\cdotp\,):S_0\to \widetilde Y$ is defined analogously.
\end{definition}

For the domain of $\mathcal{R}_M(U| m)$ we could consider any of the non-trivial 
Borel-measurable subspaces $L\subset S$, such that the
realization $M(\omega)$ belongs to $L$ almost surely. Given any two such subspaces $L_1$ and $L_2,$ the value of the function $\mathcal{R}_M(U|
m)$ can be changed in $L_1\setminus L_2$,  without contradicting the property (\ref{det. CM estimate}). 
It is tempting to choose the domain to be
the intersection of all such subspaces $L\subset S,$ but unfortunately, this intersection is the set $Z$ where the realizations of
${M}$ lie only with probability zero. It appears to be hard to pick a candidate for a smallest  space $S_0$
where the {\mntext reconstructor} $\mathcal{R}_M(U| m)$
should be defined.

Because of the above difficulties, we restricted ourselves  to the case where the operator $A$ maps $A:S\to S^1$ implying that $m\mapsto \mathcal{R}_M(U| m)$ can be defined in the whole space $S$.

\bibliographystyle{amsalpha}

\vspace{1cm}

\end{document}